\documentclass[12pt,letterpaper]{article}
\usepackage{amsthm,amssymb}
\usepackage{amsfonts}
\makeatletter
\usepackage[dvips]{color}
\usepackage{colordvi,multicol}

\newtheorem{thm}{\textbf Theorem}[section]
\newtheorem{lem}[thm]{\textbf Lemma}
\newtheorem{rem}[thm]{\textbf Remark}
\newtheorem{cor}[thm]{\textbf Corollary}

\newtheorem{prop}[thm]{\textbf Proposition}
\newtheorem{defin}[thm]{\textbf Definition}

\newtheorem{exa}[thm]{\textbf Example}
\newcommand{\be}{\begin{eqnarray*}}
\newcommand{\ee}{\end{eqnarray*}}

\begin{document}

\centerline{\Large\bf Further study on  Hunt's hypothesis (H)}
\vskip 0.2cm
\centerline{\Large\bf  for L\'{e}vy processes}

\vskip 1.6cm \centerline{Ze-Chun Hu} \centerline{\small College of Mathematics, Sichuan University, Chengdu, 610064, China}
\centerline{\small E-mail: zchu@scu.edu.cn}
\vskip 0.7cm \centerline{Wei Sun} \centerline{\small Department of
Mathematics and Statistics, Concordia University,}
\centerline{\small Montreal, H3G 1M8, Canada} \centerline{\small
E-mail: wei.sun@concordia.ca}
 \vskip 1.4cm

\vskip 0.5cm \noindent{\bf Abstract}\quad Getoor's conjecture that essentially all L\'{e}vy processes satisfy (H) is a long-standing open problem in
potential theory. In the beginning of the paper, we summarize the
main results obtained so far for the problem. Then, we present two
new necessary and sufficient conditions for the validity of (H).
Furthermore, we give applications of these new criteria. First, we
give explicit constructions of L\'evy processes satisfying (H) in
a context where previously known results could not be applied.
Second, we show that a large class of pure jump subordinators can
be decomposed into the summation of two independent subordinators
such that both of them satisfy (H).
 \smallskip

\noindent {\bf Keywords}\quad Hunt's hypothesis (H), Getoor's
conjecture, L\'{e}vy process, subordinator.

\smallskip

\noindent {\bf Mathematics Subject Classification (2010)}\quad
Primary: 60J45; Secondary: 60G51

%%%%%%%%%%%%%%%%%%%%%%%%%%%%%%%%%%%%%%%%%%%%%%%%%%%%%%%%%%%%%%%%%%%%%%%
%\tableofcontents

\section{Introduction}

A L\'evy process on $\mathbf{R}^n$ is said to satisfy Hunt's
hypothesis (H) if every semipolar set is polar. About fifty years
ago, Professor R.K. Getoor raised the problem that for which
L\'evy processes semipolar sets are always polar. His conjecture
that essentially all L\'evy processes satisfy (H) is the major
open problem in the field of potential theory for L\'evy processes
(cf. e.g. \cite[page 70]{B}).

Let us start with a brief introduction to Hunt's hypothesis (H).
For simplicity, we consider here only (H) for L\'evy processes;
however, we should point out that (H) plays a crucial role in the
potential theory of (dual) Markov processes. We refer the readers
to \cite[Chapter VI]{BG68} for a systematic introduction to (H)
for Markov processes.

Let $(\Omega,{\cal F},P)$ be a probability space and
$X=(X_t)_{t\ge 0}$ be a L\'{e}vy process on $\mathbf{R}^n$  with
L\'{e}vy-Khintchine exponent $\psi$, i.e.,
\begin{eqnarray*}
E[\exp\{i\langle z,X_t\rangle\}]=\exp\{-t\psi(z)\},\  z\in
\mathbf{R}^n,t\ge 0.
\end{eqnarray*}
Hereafter $E$ denotes the expectation {with respect to} $P$,
$\langle\cdot,\cdot\rangle$ and $|\cdot|$ denote respectively the
Euclidean inner product and norm of $\mathbf{R}^n$. The classical
L\'{e}vy-Khintchine formula tells us that
\begin{eqnarray*}
\psi(z)=i\langle a,z\rangle+\frac{1}{2}\langle
z,Qz\rangle+\int_{\mathbf{R}^n} \left(1-e^{i\langle
z,x\rangle}+i\langle z,x\rangle 1_{\{|x|<1\}}\right)\mu(dx),
\end{eqnarray*}
where $a\in \mathbf{R}^n,Q$ is a symmetric nonnegative definite
$n\times n$ matrix, and $\mu$ is a measure (called the L\'evy
measure) on $\mathbf{R}^n\backslash\{0\}$ satisfying
$\int_{\mathbf{R}^n\backslash\{0\}} (1\wedge
|x|^2)\mu(dx)<\infty$. We use Re$(\psi)$ and Im$(\psi)$ to denote
respectively the real and imaginary parts of $\psi$, and use also
$(a,Q,\mu)$ to denote $\psi$. For $x\in \mathbf{R}^n$, we denote
by $P^x$ the law of $x+X$ under $P$. In particular, $P^0=P$.

Denote by ${\cal B}^n$ the family of all nearly Borel sets of
$\mathbf{R}^n$ (cf. \cite[Definition I.10.21]{BG68}). For
$D\subset E$, we denote the first hitting time of $D$ by
$$
\sigma_D:=\inf\{t>0:X_t\in D\}.
$$
A set $D\subset E$ is called polar  if there exists a set $C\in
{\cal B}^n$ such that $D\subset C$ and $P^x(\sigma_C<\infty)=0$
for every $x\in \mathbf{R}^n$. $D$ is called a thin set if there
exists a set $C\in {\cal B}^n$ such that $D\subset C$ and
{$P^x(\sigma_C=0)=0$} for every $x\in \mathbf{R}^n$. $D$ is called
semipolar if $D\subset\bigcup_{n=1}^{\infty}D_n$ for some thin
sets $\{D_n\}_{n=1}^{\infty}$. $X$ is said to satisfy Hunt's
hypothesis (H) if every semipolar set is polar.

To appreciate the importance of (H), we recall below some
important principles of potential theory that are equivalent to
(H). We refer the readers to \cite[Proposition 1.1]{HSZ15} for a
summary of the proofs. For $\alpha>0$, a finite $\alpha$-excessive
function $f$ on $\mathbf{R}^n$ is called a regular potential
provided that $E^x\{e^{-\alpha T_n}f(X_{T_n})\}\rightarrow
E^x\{e^{-\alpha T}f(X_{T})\}$  for $x\in E$ whenever $\{T_n\}$ is
an increasing sequence of stopping times with limit $T$. Denote by
$(U^{\alpha})_{\alpha>0}$ the resolvent operators for $X$. For
each $\alpha>0$, (H) is equivalent to all of the following
principles.

\begin{itemize}
\item {\bf Bounded maximum principle:} If $\nu$ is a finite
measure with compact support $K$ such that $U^{\alpha}\nu$ is
bounded, then $\sup\{U^{\alpha}\nu(x):x\in
E\}=\sup\{U^{\alpha}\nu(x):x\in K\}$.

\item {\bf Bounded energy principle:} If $\nu$ is a finite measure
with compact support such that  $U^{\alpha}\nu$ is bounded, then
$\nu$ does not charge semipolar sets.

\item {\bf Bounded regularity principle:} If $\nu$ is a finite
measure with compact support such that $U^{\alpha}\nu$ is bounded,
then $U^{\alpha}\nu$ is regular.

\item {\bf Bounded positivity principle:} If $\nu$ is a finite
signed measure such that  $U^{\alpha}\nu$ is bounded, then $\nu
U^{\alpha}\nu\geq 0$, where $\nu
U^{\alpha}\nu:=\int_EU^{\alpha}\nu(x)\nu(dx)$.

\end{itemize}

Hunt's hypothesis (H) is also equivalent to some other important
properties of Markov processes. For example, Blumenthal and Getoor
\cite[Proposition (4.1)]{BG70} and Glover \cite[Theorem
(2.2)]{G83} showed that (H) holds if and only if the fine and
cofine topologies differ by polar sets; Fitzsimmons and Kanda
\cite{FK} showed that (H) is equivalent to the dichotomy of
capacity.

In spite of its importance, Hunt's hypothesis (H) has been
verified only in special situations. Blumenthal and Getoor
\cite{BG70} showed that all  stable processes with index
$\alpha\in (0,2)$ on the line satisfy (H). Kanda \cite{Ka76} and
Forst \cite{F75} proved independently that (H) holds if $X$ has
bounded continuous transition densities (with respect to the
Lebesgue measure $dx$) and the L\'{e}vy-Khintchine exponent $\psi$
satisfies \begin{equation}\label{Kanda}|\mbox{Im} (\psi)|\leq
M(1+\mbox{Re}(\psi)) \end{equation}
 for some constant
$M>0$. Rao \cite{R77} gave a short proof of the Kanda-Forst
theorem under the weaker condition that $X$ has resolvent
densities.
 In particular, for $n\ge 1$, all
stable processes with index $\alpha\neq 1$ satisfy (H). Kanda
\cite{Ka78} proved that (H) holds for stable processes on
$\mathbf{R}^n$ with index $\alpha= 1$ if we assume that the linear
term vanishes.  Glover and Rao \cite{GR86} proved that
$\alpha$-subordinates of general Hunt processes satisfy (H) (cf.
Proposition \ref{GR} below).
%Rao  \cite{R87}(1987) derived some of the above-mentional
%results, using some general energy considerations developed in
%Pop-Stojanovic and Rao \cite{PR81}(1981). By using the results in
%\cite{R87},
{Rao \cite{R88} proved that if all 1-excessive functions of $X$
are lower semicontinuous and \begin{equation}\label{Rao}|{\rm
Im}(\psi)|\leq (1+{\rm Re}(\psi))f(1+{\rm Re}(\psi)),
\end{equation} where $f$ is an increasing function on $[1,\infty)$
such that $\int_N^{\infty}(\lambda
f(\lambda))^{-1}d\lambda=\infty$ for every $N\geq 1$, then $X$
satisfies (H).}

In \cite{HS12}, we showed that i) if $Q$ is non-degenerate, then
$X$ satisfies (H); ii) if $Q$ is degenerate and
$\mu({\mathbf{R}^n\backslash \sqrt{Q}\mathbf{R}^n})<\infty$, then
$X$ satisfies (H) if and only if the equation
$$
\sqrt{Q}y=-a-\int_{\mathbf{R}^n\backslash
\sqrt{Q}\mathbf{R}^n}x1_{\{|x|<1\}}\mu(dx)
$$
has at least one solution $y\in \mathbf{R}^n$. In \cite{HS12}, we
also showed that if $X$ is a subordinator satisfying (H) then its
drift coefficient must be 0, i.e., $X$ must be a pure jump
subordinator. Recently, Hu, Sun and Zhang \cite{HSZ15} gave a
comparison result on L\'evy processes which implies that big jumps
have no effect on the validity of (H) in some sense. Moreover, Hu, Sun and Zhang
obtained an extended
Kanda-Forst-Rao theorem by virtue of a new necessary and sufficient condition for (H) (see
\cite[Theorems 4.3 and 4.5]{HSZ15}).

In this paper, we will continue to explore Hunt's hypothesis (H)
for L\'evy processes. In the next section, we present two new
necessary and sufficient conditions for (H), see Theorems
\ref{thm2.3} and \ref{thm2.4} below. Applications of these new
criteria for (H) will be given in Sections 3 and 4. In Section 3,
we give explicit constructions of L\'evy processes satisfying (H)
in a context where previously known results could not be applied.
In Section 4, we prove an affirmative result concerning (H) for
subordinators. We show that a large class of pure jump
subordinators can be decomposed into the summation of two
independent subordinators such that both of them satisfy (H).

%%%%%%%%%%%%%%%%%%%%%%%%%%%%%%%%%%%%%%%%%%%%%%%%%%%%%%%%%%%%%%%%%%%%%%%
\section{New necessary and sufficient conditions for (H)}\setcounter{equation}{0}

\subsection{Statement of the main theorems}

Let $X$ be a L\'evy process on $\mathbf{R}^n$ with
L\'{e}vy-Khintchine exponent $\psi$. From now on till the end of
this paper we assume that all 1-excessive functions are lower
semicontinuous, equivalently, $X$ has resolvent densities. We
refer the readers to \cite[Theorem 2.1]{Haw} for more
characterizations of this assumption.

Define
\begin{equation}\label{ABAB}
A:=1+{\rm Re}(\psi),\ \ B:=|1+\psi|.
\end{equation}
For a finite (positive) measure $\nu$ on $\mathbf{R}^n$, we denote
$$ \hat\nu(z):=\int_{\mathbf{R}^n}e^{i\langle z,x\rangle}\nu(dx).
$$
$\nu$ is said to have finite 1-energy if
\begin{equation}\label{energy}
\int_{\mathbf{R}^n}\frac{A(z)}{B^{2}(z)}|\hat \nu(z)|^2dz<\infty.
\end{equation}
Let $\nu$ be a finite measure  on $\mathbf{R}^n$ of finite
1-energy. For $\lambda>0$, we define
$$
c(\lambda):=\int_{\mathbf{R}^n}\frac{\lambda}{\lambda^2+B^2(z)}|\hat\nu(z)|^2dz.
$$
Note that
\begin{equation}\label{add98}
c(\lambda)\le\lambda\int_{\mathbf{R}^n}\frac{A(z)}{B^{2}(z)}|\hat \nu(z)|^2dz<\infty.
\end{equation}
By \cite[Theorem 1 and the proof of Theorem 2]{R88}, we have
\begin{equation}\label{ld23}
\lim_{\lambda\rightarrow\infty}c(\lambda)=\lim_{\lambda\rightarrow\infty}\int_{\mathbf{R}^n}\frac{\lambda}{\lambda^2+B^2(z)}|\hat\nu(z)|^21_{\{A(z)\le\lambda\}}dz<\infty.
\end{equation}

The remarkable result of Rao (\cite[Theorem 1]{R88}) essentially tells us that
whether $X$ satisfies (H) depends on if the limit in (\ref{ld23})
equals 0. In fact, by virtue of \cite[Theorem 1]{R88}, we proved in \cite{HSZ15} the following result.

\begin{prop}\label{thm2.1} (\cite[Theorem 5.1]{HSZ15})
(H) holds if and only if
\begin{eqnarray*}
\lim_{\lambda\rightarrow\infty}\int_{\mathbf{R}^ n}|\hat\nu(z)|^2(\lambda+{\rm
Re}\psi(z))|\lambda+\psi(z)|^{-2}dz=0
\end{eqnarray*}
for any finite measure $\nu$ of finite 1-energy.
\end{prop}

Further, by Proposition \ref{thm2.1} and \cite[the proof of
Theorem 2]{R88}, we can obtain the following result.

\begin{prop}\label{thm2.2}
(H) holds if and only if
\begin{eqnarray}\label{adds11}
\lim_{\lambda\to\infty}\int_{\mathbf{R}^n}\frac{\lambda}{\lambda^2+B^2(z)}|\hat\nu(z)|^2dz=0
\end{eqnarray}
for any finite measure $\nu$ of finite 1-energy.
\end{prop}

Let $\varsigma>1$ be a constant. We define
$$N^{\varsigma}_x:=\varsigma^{(\varsigma^x)}\ \ {\rm for}\
x\ge 0.$$ To simplify notations, we use $N$ to denote
$N^{\varsigma}$ whenever there is no confusion caused. Throughout
this paper, we use $\log$ to denote $\log_e$.

Now we state the main theorems of this section. \vskip 0.5cm
 \begin{thm}\label{thm2.3} %{\bf (New criterion for (H) of L\'evy processes)}

(i) $X$ satisfies (H)  if the following condition holds:
\begin{eqnarray}\label{zeta}
& &{ Condition\ (C^{\log})}:\ { For\ any\ finite\
measure}\ \nu\ { on}\ \mathbf{R}^n\ { of\ finite\ 1{\textrm{-}}energy,\ there}\nonumber\\
& &{exist\ a\ constant}\ \varsigma>1\  and\ a\ sequence\
\ \{y_k\uparrow\infty\}\ { such\ that}\ y_1>1\ and\nonumber\\
& &\ \ \ \ \ \ \ \ \ \ \ \ \ \ \sum_{k=1}^{\infty}\int_{\{y_k\le
B(z)< (y_k)^{\varsigma}\}} \frac{1}{B(z)\log(B(z))}|\hat
\nu(z)|^2dz<\infty.
\end{eqnarray}

(ii) Suppose $X$ satisfies (H). Then, for\ any\ finite\ measure
$\nu$ on $\mathbf{R}^n$  of\ finite\ 1-energy and any
$\varsigma>1$, there exists a sequence\ $\{y_k\uparrow\infty\}$
such\ that $y_1>1$ and (\ref{zeta})\ {holds}.

\end{thm}
\vskip 0.5cm

\begin{thm}\label{thm2.4}
 (i)
$X$ satisfies (H)  if the following condition holds:
\begin{eqnarray}\label{new345}
& &{ Condition\ (C^{\log\log})}:\ { For\ any\ finite\
measure}\ \nu\ { on}\ \mathbf{R}^n\ { of\ finite\ 1{\textrm{-}}energy,}\nonumber\\
& &\ \ \ \ \ {there\ exist\ a\ constant}\ \varsigma>1\ { and\ a\
sequence\ of\ positive\ numbers}\
\{x_k\}\nonumber\\
& &\ \ \ \ \ { such\ that}\ N^{\varsigma}_{x_1}>e,\
x_k+1<x_{k+1},\ k\in
\mathbf{N},\ \sum_{k=1}^{\infty}\frac{1}{x_k}=\infty,\ { and}\nonumber\\
& &\ \ \ \ \ \ \sum_{k=1}^{\infty}\int_{\{N^\varsigma_{x_k}\le
B(z)< N^\varsigma_{x_k+1}\}} \frac{1}{B(z)\log
(B(z))[\log\log(B(z))]}|\hat \nu(z)|^2dz<\infty.
\end{eqnarray}

(ii) Suppose $X$ satisfies (H). Then, for any finite measure $\nu$
on $\mathbf{R}^n$ of finite 1-energy and any $\varsigma>1$, there
exists a sequence of positive numbers $\{x_k\}$ such that $
N^{\varsigma}_{x_1}>e$, $x_k+1<x_{k+1}$, $k\in \mathbf{N}$, $
\sum_{k=1}^{\infty}\frac{1}{x_k}=\infty$, and (\ref{new345})\
{holds}.

\end{thm}

\subsection{Proof of Theorem \ref{thm2.3}}
Before proving the main theorems, we first prove a lemma which has independent interest.  \vskip 0.5cm

\begin{lem}\label{slight2}
Let $\nu$ be a finite measure  on $\mathbf{R}^n$ of finite 1-energy. Then, the following condition
is fulfilled for any $\delta>0$.
\begin{eqnarray}\label{slight2-0}
{Condition\ (C^{\delta})}:\ \int_{\mathbf{R}^n}
\frac{1}{B(z)\log(2+B(z))[\log\log(2+B(z))]^{1+\delta}}|\hat
\nu(z)|^2dz<\infty.\ \ \ \
\end{eqnarray}
\end{lem}

\noindent \emph{Proof.} We fix a $\delta>0$ and let $\nu$ be an
arbitrary finite measure on $\mathbf{R}^n$ of finite 1-energy.
Denote $$ F:=\{z\in\mathbf{R}^n:B(z)\ge 4A(z)\}.
$$
Then, we obtain by the fact that $B(z)\geq 1$ and (\ref{energy}) that
\begin{eqnarray}\label{slight2-a}
&&\int_{\mathbf{R}^n\backslash F}\frac{1}{B(z)\log(2+B(z))[\log\log(2+B(z))]^{1+\delta}}|\hat
\nu(z)|^2dz\nonumber\\
&\leq& \frac{4}{\log3(\log\log3)^{1+\delta}}\int_{\mathbf{R}^n\backslash F}\frac{A(z)}{B^2(z)}|\hat
\nu(z)|^2dz\nonumber\\
&<&\infty.
\end{eqnarray}

By (\ref{slight2-a}), to prove (\ref{slight2-0}), it suffices to show that
$$
\int_{F} \frac{1}{B(z)\log(2+B(z))[\log\log(2+B(z))]^{1+\delta}}|\hat
\nu(z)|^2dz<\infty.$$
Further, it suffices to show that
\begin{eqnarray}\label{slight2-a-1}
\int_{F} \frac{1}{B(z)\log(B(z))[\log\log(B(z))]^{1+\delta}}|\hat
\nu(z)|^2dz<\infty.
\end{eqnarray}

Set $N_k=2^{(2^k)}$ for $k\ge 1$. Then, for each $k\geq 2$, we have
\begin{eqnarray}\label{slight2-c}
\log(N_{k-1})&=&\log(N_k)-\log(N_{k-1}).
\end{eqnarray}
By the fact that $\log(2/\log2)<4(\log\log4)=4\log\log(N_1)$, we get
\begin{eqnarray}
k&=&\log_2\log_2(N_{k-1})^2\nonumber\\
&=&\frac{1}{\log2}\left(\log\log( N_{k-1})+\log\left(\frac{2}{\log2}\right)\right)\nonumber\\
&<& \frac{5}{\log2}\log\log(N_{k-1}).\label{slight2-d}
\end{eqnarray}

Define $f^{\delta}(\lambda)=1$
when $1\le\lambda<4$ and
$$
f^{\delta}(\lambda)=k^{1+\delta}(\log(N_k)-\log(N_{k-1})),\ \ {\rm
when}\ N_{k-1}\le \lambda<N_k,\ k\ge 2.
$$
Then, we have
\begin{eqnarray}\label{slight2-b}
\int_1^{\infty}\frac{d\lambda}{\lambda
f^{\delta}(\lambda)}
&=&\int_1^4\frac{1}{\lambda}d\lambda+
\sum_{k=2}^{\infty}\int_{N_{k-1}}^{N_k}\frac{1}{k^{1+\delta}(\log(N_k)-\log(N_{k-1}))\lambda}d\lambda\nonumber\\
&=&\log4+\sum_{k=2}^{\infty}\frac{1}{k^{1+\delta}}\nonumber\\
&<&\infty.
\end{eqnarray}

By (\ref{add98}), (\ref{ld23}), (\ref{slight2-b}),  Fubini's theorem, the definition of the function $f^{\delta}(\lambda)$, (\ref{slight2-c}) and (\ref{slight2-d}),    we obtain that
\begin{eqnarray*}\label{wei22}
\infty&>&\int_1^{\infty}\frac{d\lambda}{\lambda
f^{\delta}(\lambda)}\int_{\mathbf{R}^n}\frac{\lambda}{\lambda^2+B^2(z)}|\hat\nu(z)|^21_{\{A(z)\le\lambda\}}dz\nonumber\\
&=&\int_{\mathbf{R}^n}|\hat\nu(z)|^2dz\int_{A(z)}^{\infty}\frac{d\lambda}{f^{\delta}(\lambda)
(\lambda^2+B^2(z))}\nonumber\\
&=&\int_{\mathbf{R}^n}|\hat\nu(z)|^2dz\int_{1}^{\infty}\frac{A(z)d\eta}{f^{\delta}(A(z)\eta)
(A^2(z)\eta^2+B^2(z))}\quad\quad ({\rm by\ letting}\ \lambda=A(z)\eta)\nonumber\\
&=&\int_{\mathbf{R}^n}\frac{A(z)}{B^{2}(z)}|\hat\nu(z)|^2dz\int_{1}^{\infty}\frac{(\frac{B(z)}{A(z)})^2d\eta}{f^{\delta}(A(z)\eta)
(\eta^2+(\frac{B(z)}{A(z)})^2)}\nonumber\\
&\ge&\int_{\mathbf{R}^n}\frac{A(z)}{B^{2}(z)}|\hat\nu(z)|^2dz\sum_{k=2}^{\infty}\int_{\{1\vee\frac{N_{k-1}}{A(z)}\le
\eta<
\frac{N_k}{A(z)}\}}\frac{(\frac{B(z)}{A(z)})^2d\eta}{k^{1+\delta}(\log(N_k)-\log(N_{k-1}))
(\eta^2+(\frac{B(z)}{A(z)})^2)}\nonumber\\
&=&\int_{\mathbf{R}^n}\frac{1}{B(z)}|\hat\nu(z)|^2dz
\sum_{k=2}^{\infty}\frac{1}{k^{1+\delta}(\log(N_k)-\log(N_{k-1}))}
\int_{\{1\vee\frac{N_{k-1}}{A(z)}\le
\eta<
\frac{N_k}{A(z)}\}}\frac{\frac{A(z)}{B(z)}d\eta}{
1+(\frac{A(z)}{B(z)}\eta)^2}\nonumber\\
&=&\int_{\mathbf{R}^n}\frac{1}{B(z)}|\hat\nu(z)|^2\sum_{k=2}^{\infty}\frac{\arctan\frac{N_k}{B(z)}-
\arctan\frac{N_{k-1}\vee A(z)}{B(z)}}{k^{1+\delta}(\log(N_k)-\log(N_{k-1}))}dz\nonumber\\
&\ge&\sum_{k_0=2}^{\infty}\int_{\{N_{k_0-1}\le B(z)<
N_{k_0}\}}\frac{1}{B(z)}|\hat\nu(z)|^2\sum_{k=2}^{\infty}\frac{\arctan\frac{N_k}{B(z)}-
\arctan\frac{N_{k-1}\vee A(z)}{B(z)}}{k^{1+\delta}(\log(N_k)-\log(N_{k-1}))}dz\nonumber\\
&\ge&\sum_{k_0=2}^{\infty}\int_{\{z\in F:\,N_{k_0-1}\le B(z)<
N_{k_0}\}}\frac{1}{B(z)}|\hat\nu(z)|^2\frac{\arctan\frac{N_{k_0}}{B(z)}-
\arctan(\frac{N_{k_0-1}}{B(z)}\vee\frac{1}{4})}{k_0^{1+\delta}(\log(N_{k_0})-\log(N_{k_0-1}))}dz\nonumber\\
&=&\sum_{k_0=2}^{\infty}\left\{\int_{\{z\in F:\,N_{k_0-1}\le B(z)<
\frac{N_{k_0}}{2}\}}\frac{1}{B(z)}|\hat\nu(z)|^2\frac{\arctan\frac{N_{k_0}}{B(z)}-
\arctan(\frac{N_{k_0-1}}{B(z)}\vee\frac{1}{4})}{k_0^{1+\delta}(\log(N_{k_0})-\log(N_{k_0-1}))}dz\right.\nonumber\\
&&+\left.\int_{\{z\in F:\,\frac{N_{k_0}}{2}\le B(z)<
N_{k_0}\}}\frac{1}{B(z)}|\hat\nu(z)|^2\frac{\arctan\frac{N_{k_0}}{B(z)}-
\arctan(\frac{N_{k_0-1}}{B(z)}\vee\frac{1}{4})}{k_0^{1+\delta}(\log(N_{k_0})-\log(N_{k_0-1}))}dz\right\}\nonumber\\
&\ge&\sum_{k_0=2}^{\infty}\left\{\int_{\{z\in F:\,N_{k_0-1}\le
B(z)<
\frac{N_{k_0}}{2}\}}\frac{1}{B(z)}|\hat\nu(z)|^2\frac{\arctan2-
\arctan1}{k_0^{1+\delta}(\log(N_{k_0})-\log(N_{k_0-1}))}dz\right.\nonumber\\
&&+\left.\int_{\{z\in F:\,\frac{N_{k_0}}{2}\le B(z)<
N_{k_0}\}}\frac{1}{B(z)}|\hat\nu(z)|^2\frac{\arctan1-
\arctan\frac{1}{2}}{k_0^{1+\delta}(\log(N_{k_0})-\log(N_{k_0-1}))}dz\right\}\nonumber\\
&\ge&\sum_{k_0=2}^{\infty}\int_{\{z\in F:\,N_{k_0-1}\le B(z)<
{N_{k_0}}\}}\frac{1}{B(z)}|\hat\nu(z)|^2\frac{\arctan1-
\arctan\frac{1}{2}}{k_0^{1+\delta}\log(N_{k_0-1})}dz\nonumber\\
&\ge&\sum_{k_0=2}^{\infty}\int_{\{z\in F:\,N_{k_0-1}\le B(z)<
{N_{k_0}}\}}\frac{1}{B(z)}|\hat\nu(z)|^2\frac{\arctan1-
\arctan\frac{1}{2}}{(\frac{5}{\log2}\log\log (N_{k_0-1}))^{1+\delta}\log(N_{k_0-1})}dz\nonumber\\
&\ge&\frac{\arctan1-
\arctan\frac{1}{2}}{(5/\log2)^{1+\delta}}\sum_{k_0=2}^{\infty}\int_{\{z\in
F:\,N_{k_0-1}\le B(z)<
{N_{k_0}}\}}\frac{1}{B(z)\log(B(z))[\log\log(B(z))]^{1+\delta}}|\hat\nu(z)|^2dz\nonumber\\
&=&\frac{\arctan1-
\arctan\frac{1}{2}}{(5/\log2)^{1+\delta}}\int_{F}\frac{1}{B(z)\log(B(z))
[\log\log(B(z))]^{1+\delta}}|\hat\nu(z)|^2dz.
\end{eqnarray*}
Therefore,  (\ref{slight2-a-1}) holds and the proof is complete. \hfill\fbox

\bigskip

\noindent{\bf Proof of Theorem \ref{thm2.3}.}

(i) Suppose Condition ${\rm (C^{\log})}$ holds. We will show that $X$ satisfies (H). By Proposition \ref{thm2.2}, we need prove that (\ref{adds11}) holds for any finite measure $\nu$ of finite 1-energy.

Let $\nu$ be a finite measure on $\mathbf{R}^n$ of finite
1-energy. By Condition ${\rm (C^{\log})}$, there exist a constant
$\varsigma>1$ and a sequence $\{y_k\uparrow\infty\}$ such that
$y_1>1$ and (\ref{zeta}) holds. We assume without loss of
generality that $y_1>\varsigma$ and $(y_k)^{\varsigma}<y_{k+1}$
for $k\in \mathbf{N}$. Set
$$
x_k=\log_{\varsigma}(\log_{\varsigma}{y_k}),\ \
k\in\mathbf{N}.
$$
For
each $k\in\mathbf{N}$, we have
\begin{equation}\label{sep1}
x_k+1=\log_{\varsigma}(\log_{\varsigma}{(y_k)^\varsigma})<\log_{\varsigma}(\log_{\varsigma}y_{k+1})=x_{k+1}.\end{equation}

Note that when $B(z)<N_{x_1}$ we have
\begin{eqnarray}\label{proof-them2.3-a}
\frac{\lambda}{\lambda^2+B^2(z)}\le\frac{1}{2B(z)}<
\frac{N_{x_1}}{2}\frac{A(z)}{B^2(z)}.
\end{eqnarray}
By (\ref{energy}), (\ref{proof-them2.3-a}) and the dominated convergence
theorem, we get
\begin{eqnarray*}
\lim_{\lambda\rightarrow\infty}\int_{\{B(z)<
N_{x_1}\}}\frac{\lambda}{\lambda^2+B^2(z)}|\hat\nu(z)|^2dz=0.
\end{eqnarray*}
Hence, to prove (\ref{adds11}), it suffices to prove that
\begin{eqnarray}\label{proof-them2.3-b}
\lim_{\lambda\rightarrow\infty}\int_{\{B(z)\ge
N_{x_1}\}}\frac{\lambda}{\lambda^2+B^2(z)}|\hat\nu(z)|^2dz=0.
\end{eqnarray}

For $\lambda>1$, we define
$$
g(\lambda):=\sum_{k=1}^{\infty}\frac{1_{\{N_{x_k+\frac{1}{2}}\le\lambda<
N_{x_k+\frac{3}{4}}\}}}{\log\lambda}.
$$
Then, we have
\begin{eqnarray}\label{proof-them2.3-c}
\int_{N_{x_1}}^{\infty}\frac{g(\lambda)}{\lambda}d\lambda
&=&\sum_{k=1}^{\infty}\int_{N_{x_k+\frac{1}{2}}}^{N_{x_k+\frac{3}{4}}}\frac{1}{\lambda\log \lambda}d\lambda\nonumber\\
&=&\sum_{k=1}^{\infty}(\log\log(N_{x_k+\frac{3}{4}})-\log\log(N_{x_k+\frac{1}{2}}))\nonumber\\
&=&\sum_{k=1}^{\infty}\left(\log\log(\varsigma^{(\varsigma^{x_k+\frac{3}{4}})}\right)-
\log\log\left(\varsigma^{(\varsigma^{x_k+\frac{1}{2}})})\right)\nonumber\\
&=&\sum_{k=1}^{\infty}\frac{1}{4}\log\varsigma\nonumber\\
&=&\infty.
\end{eqnarray}

By Fubini's theorem, the definition of the function $g(\lambda)$, and the inequality that $\frac{\pi}{2}-\arctan x\le\frac{1}{x}$
for $x>0$,  we obtain that
\begin{eqnarray}\label{proof-them2.3-d}
& &\int_{N_{x_1}}^{\infty}\frac{g(\lambda)d\lambda}{\lambda
}\int_{\{B(z)\ge
N_{x_1}\}}\frac{\lambda}{\lambda^2+B^2(z)}|\hat\nu(z)|^2dz\nonumber\\
&=&\int_{\{B(z)\ge
N_{x_1}\}}|\hat\nu(z)|^2dz\int_{N_{x_1}}^{\infty}\frac{g(\lambda)d\lambda}{
\lambda^2+B^2(z)}\nonumber\\
&=&\int_{\{B(z)\ge N_{x_1}\}}|\hat\nu(z)|^2dz\int_{
\frac{N_{x_1}}{A(z)}}^{\infty}\frac{g(A(z)\eta)A(z)d\eta}{
A^2(z)\eta^2+B^2(z)}\quad\quad ({\rm by\ letting}\ \lambda=A(z)\eta)\nonumber\\
&=&\int_{\{B(z)\ge
N_{x_1}\}}\frac{A(z)}{B^{2}(z)}|\hat\nu(z)|^2dz\int_{
\frac{N_{x_1}}{A(z)}}^{\infty}\frac{g(A(z)\eta)(\frac{B(z)}{A(z)})^2d\eta}{
\eta^2+(\frac{B(z)}{A(z)})^2}\nonumber\\
&\le&\int_{\{B(z)\ge
N_{x_1}\}}\frac{A(z)}{B^{2}(z)}|\hat\nu(z)|^2dz\sum_{k=1}^{\infty}\int_{\{\frac{N_{x_k+\frac{1}{2}}}{A(z)}\le
\eta<
\frac{N_{x_k+\frac{3}{4}}}{A(z)}\}}\frac{(\frac{B(z)}{A(z)})^2d\eta}{\log(N_{x_k+\frac{1}{2}})
(\eta^2+(\frac{B(z)}{A(z)})^2)}\nonumber\\
&=&\int_{\{B(z)\ge
N_{x_1}\}}\frac{1}{B(z)}|\hat\nu(z)|^2dz
\sum_{k=1}^{\infty}\frac{1}{\log(N_{x_k+\frac{1}{2}})}\int_{\{\frac{N_{x_k+\frac{1}{2}}}{A(z)}\le
\eta<
\frac{N_{x_k+\frac{3}{4}}}{A(z)}\}}\frac{\frac{A(z)}{B(z)}d\eta}{1+(\frac{A(z)}{B(z)}\eta)^2}\nonumber\\
&=&\int_{\{B(z)\ge
N_{x_1}\}}\frac{1}{B(z)}|\hat\nu(z)|^2\sum_{k=1}^{\infty}\frac{\arctan\frac{N_{x_k+\frac{3}{4}}}{B(z)}-
\arctan\frac{N_{x_k+\frac{1}{2}}}{B(z)}}{\log(N_{x_k+\frac{1}{2}})}dz\nonumber\\
&=&\sum_{l=1}^{\infty}\int_{\{N_{x_{l}}\le B(z)<
N_{(x_{l})+1}\}}\frac{1}{B(z)}|\hat\nu(z)|^2\sum_{k=1}^{\infty}\frac{\arctan\frac{N_{x_k+\frac{3}{4}}}{B(z)}-
\arctan\frac{N_{x_k+\frac{1}{2}}}{B(z)}}{\log(N_{x_k+\frac{1}{2}})}dz\nonumber\\
& &+\sum_{l=1}^{\infty}\int_{\{N_{(x_l)+1}\le B(z)<
N_{x_{(l+1)}}\}}\frac{1}{B(z)}|\hat\nu(z)|^2\sum_{k=1}^{\infty}\frac{\arctan\frac{N_{x_k+\frac{3}{4}}}{B(z)}-
\arctan\frac{N_{x_k+\frac{1}{2}}}{B(z)}}{\log(N_{x_k+\frac{1}{2}})}dz\nonumber\\
  &=&\sum_{l=1}^{\infty}\int_{\{N_{x_{l}}\le B(z)<
N_{(x_{l})+1}\}}\frac{1}{B(z)}|\hat\nu(z)|^2
\left\{\sum_{k=1}^{l-1}\frac{\arctan\frac{N_{x_k+\frac{3}{4}}}{B(z)}-
\arctan\frac{N_{x_k+\frac{1}{2}}}{B(z)}}{\log(N_{x_k+\frac{1}{2}})}\right.\nonumber\\
& &\ \ \ \
\left.+\frac{\arctan\frac{N_{x_l+\frac{3}{4}}}{B(z)}-
\arctan\frac{N_{x_l+\frac{1}{2}}}{B(z)}}{\log(N_{x_l+\frac{1}{2}})}+\sum_{k=l+1}^{\infty}\frac{\arctan\frac{N_{x_k+\frac{3}{4}}}{B(z)}-
\arctan\frac{N_{x_k+\frac{1}{2}}}{B(z)}}{\log(N_{x_k+\frac{1}{2}})}\right\}dz\nonumber\\
&&+\sum_{l=1}^{\infty}\int_{\{N_{(x_l)+1}\le B(z)<
N_{x_{(l+1)}}\}}\frac{1}{B(z)}|\hat\nu(z)|^2
\left\{\sum_{k=1}^{l}\frac{\arctan\frac{N_{x_k+\frac{3}{4}}}{B(z)}-
\arctan\frac{N_{x_k+\frac{1}{2}}}{B(z)}}{\log(N_{x_k+\frac{1}{2}})}\right.\nonumber\\
& &\ \ \ \
\left.+\sum_{k=l+1}^{\infty}\frac{\arctan\frac{N_{x_k+\frac{3}{4}}}{B(z)}-
\arctan\frac{N_{x_k+\frac{1}{2}}}{B(z)}}{\log(N_{x_k+\frac{1}{2}})}\right\}dz\nonumber\\
&\le&\sum_{l=1}^{\infty}\int_{\{N_{x_{l}}\le B(z)<
N_{(x_{l})+1}\}}\frac{1}{B(z)}|\hat\nu(z)|^2\left\{\frac{N_{x_{(l-1)}+\frac{3}{4}}}{B(z)}
\sum_{k=1}^{l-1}\frac{1}{\log(N_{x_k+\frac{1}{2}})}\right.\nonumber\\
& &\ \ \ \
\left.+\frac{\pi}{2\log(N_{x_l+\frac{1}{2}})}+\sum_{k=l+1}^{\infty}\frac{\frac{B(z)}{N_{x_k+\frac{1}{2}}}}
{\log(N_{x_k+\frac{1}{2}})}\right\}dz\nonumber\\
&&+\sum_{l=1}^{\infty}\int_{\{N_{(x_l)+1}\le B(z)<
N_{x_{(l+1)}}\}}\frac{1}{B(z)}|\hat\nu(z)|^2\left\{\frac{N_{(x_{l})+\frac{3}{4}}}{B(z)}
\sum_{k=1}^{l}\frac{1}{\log(N_{x_k+\frac{1}{2}})}\right.\nonumber\\
& &\ \ \ \
\left.+\sum_{k=l+1}^{\infty}\frac{\frac{B(z)}{N_{x_k+\frac{1}{2}}}}
{\log(N_{x_k+\frac{1}{2}})}\right\}dz.
\end{eqnarray}

By (\ref{sep1}), we get $
x_k+\frac{1}{2}\geq (x_1+\frac{1}{2})+(k-1)$ for $k\geq 1$.
Hence
\begin{eqnarray}\label{proof-them2.3-e}
\sum_{k=1}^{l-1}\frac{1}{\log(N_{x_k+\frac{1}{2}})}&\le&\sum_{k=1}^{\infty}\frac{1}{\log(N_{x_k+\frac{1}{2}})}\nonumber\\
&=&
\sum_{k=1}^{\infty}\frac{1}{\log(\varsigma^{(\varsigma^{x_k+\frac{1}{2}})})}\nonumber\\
&=&\sum_{k=1}^{\infty}\frac{1}{\varsigma^{x_k+\frac{1}{2}}\log\varsigma}\nonumber\\
&=&\frac{1}{\varsigma^{x_1+\frac{1}{2}}\log\varsigma}
\sum_{k=1}^{\infty}\frac{1}{\varsigma^{k-1}}\nonumber\\
&=& \frac{1}{\varsigma^{x_1+\frac{1}{2}}(\log\varsigma)(1-\frac{1}{\varsigma})}.
\end{eqnarray}
When $N_{x_{l}}\le B(z)$, we obtain by (\ref{sep1}) and the definition of $N_x$ that
\begin{eqnarray}\label{proof-them2.3-f}
\frac{N_{x_{(l-1)}+\frac{3}{4}}}{B(z)}
&\leq& \frac{N_{x_l-\frac{1}{4}}}{B(z)}\nonumber\\
&=&\frac{(N_{x_l})^{\varsigma^{-\frac{1}{4}}}}{B(z)}\nonumber\\
&\le&\frac{1}{(B(z))^{1-\varsigma^{-\frac{1}{4}}}}.
\end{eqnarray}
By (\ref{proof-them2.3-e}) and (\ref{proof-them2.3-f}), we obtain that when $N_{x_{l}}\le B(z)$,
\begin{eqnarray}\label{proof-them2.3-g}
\frac{N_{x_{(l-1)}+\frac{3}{4}}}{B(z)}\sum_{k=1}^{l-1}\frac{1}{\log(N_{x_k+\frac{1}{2}})}
\leq \frac{1}{\varsigma^{x_1+\frac{1}{2}}(\log\varsigma)(1-\frac{1}{\varsigma})(B(z))^{1-\varsigma^{-\frac{1}{4}}}}.
\end{eqnarray}
When $N_{x_{l}}\le B(z)<N_{(x_l)+1}$, by the definition of $N_x$,
we get
\begin{eqnarray}\label{proof-them2.3-h}
\frac{\pi}{2\log(N_{x_l+\frac{1}{2}})}
&=&\frac{\pi}{2\log(\varsigma^{(\varsigma^{x_l+\frac{1}{2}})})}
\nonumber\\
&=&\frac{\pi}{2\varsigma^{-\frac{1}{2}}\log(\varsigma^{(\varsigma^{(x_l)+1})})}\nonumber\\
&=&\frac{\pi}{2\varsigma^{-\frac{1}{2}}\log(N_{(x_l)+1})}\nonumber\\
&<&\frac{\pi \varsigma^{\frac{1}{2}}}{2\log
(B(z))}.
\end{eqnarray}
When $B(z)<N_{(x_l)+1}$, by (\ref{sep1}), (\ref{proof-them2.3-e}) and the definition of $N_x$, we get
\begin{eqnarray}\label{proof-them2.3-i}
\sum_{k=l+1}^{\infty}\frac{\frac{B(z)}{N_{x_k+\frac{1}{2}}}}
{\log(N_{x_k+\frac{1}{2}})}&\leq&\frac{B(z)}{N_{x_{(l+1)}+\frac{1}{2}}}\sum_{k=1}^{\infty}\frac{1}
{\log(N_{x_k+\frac{1}{2}})}\nonumber\\
&\leq&\frac{B(z)}{N_{(x_l)+1+\frac{1}{2}}}\sum_{k=1}^{\infty}\frac{1}
{\log(N_{x_k+\frac{1}{2}})}\nonumber\\
&=&\frac{B(z)}{(N_{(x_l)+1})^{\varsigma^{\frac{1}{2}}}}\cdot
\frac{1}{\varsigma^{x_1+\frac{1}{2}}(\log\varsigma)(1-\frac{1}{\varsigma})}\nonumber\\
&\leq&\frac{1}{\varsigma^{x_1+\frac{1}{2}}(\log\varsigma)(1-\frac{1}{\varsigma})(B(z))^{\varsigma^{\frac{1}{2}}-1}}.
\end{eqnarray}
Similar to (\ref{proof-them2.3-g}) and (\ref{proof-them2.3-i}), we
can show that when $N_{(x_l)+1}\le B(z)<N_{x_{(l+1)}}$,
\begin{eqnarray}\label{proof-them2.3-j}
\frac{N_{(x_l)+\frac{3}{4}}}{B(z)}\sum_{k=1}^{l-1}\frac{1}{\log(N_{x_k+\frac{1}{2}})}
\leq \frac{1}{\varsigma^{x_1+\frac{1}{2}}(\log\varsigma)(1-\frac{1}{\varsigma})(B(z))^{1-\varsigma^{-\frac{1}{4}}}},
\end{eqnarray}
and
\begin{eqnarray}\label{proof-them2.3-k}
\sum_{k=l+1}^{\infty}\frac{\frac{B(z)}{N_{x_k+\frac{1}{2}}}}
{\log(N_{x_k+\frac{1}{2}})}
\leq\frac{1}{\varsigma^{x_1+\frac{1}{2}}(\log\varsigma)(1-\frac{1}{\varsigma})(B(z))^{\varsigma^{\frac{1}{2}}-1}}.
\end{eqnarray}

We fix a $\delta>0$. Note that
$$
\lim_{\eta\to\infty}\frac{\frac{1}{\eta^{1-\varsigma^{-\frac{1}{4}}}}+
\frac{1}{\eta^{\varsigma^{\frac{1}{2}}-1}}}{\frac{1}{\log(2+\eta)[\log\log(2+\eta)]^{1+\delta}}}=0,
$$
and $N_{x_l}=y_l$, $N_{(x_l)+1}=(y_l)^{\varsigma}$ for $l\ge 1$. We obtain by (\ref{proof-them2.3-d}), (\ref{proof-them2.3-g})-(\ref{proof-them2.3-k}), (\ref{zeta}) and Lemma \ref{slight2} that
\begin{eqnarray}\label{proof-them2.3-l}
& &\int_{N_{x_1}}^{\infty}\frac{g(\lambda)d\lambda}{\lambda
}\int_{\{B(z)\ge
N_{x_1}\}}\frac{\lambda}{\lambda^2+B^2(z)}|\hat\nu(z)|^2dz\nonumber\\
&\le&\frac{1}{\varsigma^{x_1+\frac{1}{2}}(\log\varsigma)(1-\frac{1}{\varsigma})}
\int_{\mathbf{R}^n}\frac{1}{B(z)}|\hat\nu(z)|^2
\left(\frac{1}{(B(z))^{1-\varsigma^{-\frac{1}{4}}}}+
\frac{1}{(B(z))^{\varsigma^{\frac{1}{2}}-1}}\right)dz\nonumber\\
&&+\frac{\pi
\varsigma^{\frac{1}{2}}}{2}\sum_{l=1}^{\infty}\int_{\{N_{x_{l}}\le
B(z)< N_{(x_{l})+1}\}}\frac{1}{B(z)\log(B(z))}|\hat\nu(z)|^2dz\nonumber\\
&\le&D_1\int_{\mathbf{R}^n}
\frac{1}{B(z)\log(2+B(z))[\log\log(2+B(z))]^{1+\delta}}|\hat
\nu(z)|^2dz\nonumber\\
&&+ \frac{\pi
\varsigma^{\frac{1}{2}}}{2}\sum_{l=1}^{\infty}\int_{\{y_l\le
B(z)< (y_l)^{\varsigma}\}}\frac{1}{B(z)\log(B(z))}|\hat\nu(z)|^2dz\nonumber\\
&<&\infty,
\end{eqnarray}
where $D_1$ is a positive constant depending only on
$\varsigma$.

By (\ref{proof-them2.3-c}) and (\ref{proof-them2.3-l}), we obtain
(\ref{proof-them2.3-b}). Therefore, the proof of (i) is complete.
\vskip 0.2cm
 (ii) Suppose that $X$ satisfies (H). Let $\nu$ be an
arbitrary finite measure on $\mathbf{R}^n$ of finite 1-energy and
$\varsigma>1$ be a constant. By Proposition \ref{thm2.2},
$\lim_{\lambda\rightarrow\infty}c(\lambda)=0$. We choose an
increasing sequence of positive numbers $\{x_k\}$ such that
\begin{eqnarray}\label{proof-them2.3-ii-a}
c(\lambda)\le\frac{1}{2^k},\ \ {\rm if}\ \lambda\ge x_k.
\end{eqnarray}
We assume without loss of generality that $x_k+1<x_{k+1},\
k\in\mathbf{N}$. Set
$$
y_k=N_{x_k}=\varsigma^{(\varsigma^{x_k})}\ \ {\rm for}\ k\in\mathbf{N}. $$Denote $$
F:=\{z\in\mathbf{R}^n:B(z)\ge 2A(z)\}.
$$

Since
\begin{eqnarray}\label{bn111}
&&\sum_{k=1}^{\infty}\int_{\{z\in \mathbf{R}^n\backslash F:N_{x_k}\le B(z)< N_{x_k+1}\}}
\frac{1}{B(z)\log(B(z))}|\hat \nu(z)|^2dz\nonumber\\
&=& \int_{\{z\in \mathbf{R}^n\backslash F:N_{x_1}\le B(z)\}}
\frac{1}{B(z)\log(B(z))}|\hat \nu(z)|^2dz\nonumber\\
&=&\int_{\{z\in \mathbf{R}^n\backslash F:N_{x_1}\le B(z)\}}
\frac{A(z)}{B^2(z)}\cdot \frac{B(z)}{A(z)\log(B(z))}|\hat \nu(z)|^2dz\nonumber\\
&\le& \frac{2}{\log (N_{x_1})}\int_{\mathbf{R}^n}
\frac{A(z)}{B^2(z)}|\hat \nu(z)|^2dz\nonumber\\
&<&\infty,
\end{eqnarray}
to prove (\ref{zeta}), it suffices to prove that
\begin{eqnarray}\label{proof-them2.3-ii-a-1}
\sum_{k=1}^{\infty}\int_{\{z\in F:N_{x_k}\le B(z)< N_{x_k+1}\}}
\frac{1}{B(z)\log(B(z))}|\hat \nu(z)|^2dz<\infty.
\end{eqnarray}

For $\lambda>1$, we define $ f(\lambda)=\log\lambda$. Set $$
\Lambda:=\bigcup_{k=1}^{\infty}\{\lambda:N_{x_{k}}\le
\lambda<N_{x_k+1}\}.
$$
Then, by (\ref{proof-them2.3-ii-a}), we get
\begin{eqnarray}\label{proof-them2.3-ii-b}
\int_{\Lambda}\frac{c(\lambda)}{\lambda
f(\lambda)}d\lambda&=&\sum_{k=1}^{\infty}\int_{N_{x_k}}^{N_{x_k+1}}\frac{c(\lambda)}{\lambda
\log\lambda} d\lambda\nonumber\\
&\le&\sum_{k=1}^{\infty}\frac{1}{2^k}\int_{N_{x_k}}^{N_{x_k+1}}\frac{1}{\lambda
\log\lambda} d\lambda\nonumber\\
&=&\sum_{k=1}^{\infty}\frac{1}{2^k}[\log\log(N_{x_k+1})-\log\log(N_{x_k})]\nonumber\\
&=&\sum_{k=1}^{\infty}\frac{1}{2^k}[\log\log(\varsigma^{(\varsigma^{x_k+1})})-\log\log(\varsigma^{(\varsigma^{x_k})})]\nonumber\\
&=&\log\varsigma\sum_{k=1}^{\infty}\frac{1}{2^k}\nonumber\\
&<&\infty.
\end{eqnarray}

For $k\geq 1$, when $N_{x_{k}}\le B(z)<
N_{x_{k}+1}$,  we have
\begin{eqnarray}\label{proof-them2.3-ii-c}
\frac{1}{\log(N_{x_{k}+1})}&=&\frac{1}{\log(\varsigma^{(\varsigma^{x_k+1})})}\nonumber\\
&=&\frac{1}{\varsigma\log(\varsigma^{(\varsigma^{x_k})})}\nonumber\\
&=&\frac{1}{\varsigma\log(N_{x_{k}})}\nonumber\\
&\geq& \frac{1}{\varsigma\log(B(z))}.
\end{eqnarray}

By (\ref{proof-them2.3-ii-b}), Fubini's theorem, (\ref{proof-them2.3-ii-c}), $N_{x_k}=y_k$ and $N_{(x_k)+1}=(y_k)^{\varsigma}$ for $k\geq 1$,  we obtain
\begin{eqnarray}
\infty&>&\int_{\Lambda}\frac{d\lambda}{\lambda
f(\lambda)}\int_{\mathbf{R}^n}\frac{\lambda}{\lambda^2+B^2(z)}|\hat\nu(z)|^21_{\{A(z)\le\lambda\}}dz\nonumber\\
&=&\int_{\mathbf{R}^n}|\hat\nu(z)|^2dz\int_{\Lambda\cap\{A(z)\le\lambda\}}\frac{d\lambda}{f(\lambda)
(\lambda^2+B^2(z))}\nonumber\\
&=&\int_{\mathbf{R}^n}|\hat\nu(z)|^2dz\int_{\frac{\Lambda}{A(z)}\cap\{1\le\eta\}}\frac{A(z)d\eta}{f(A(z)\eta)
(A^2(z)\eta^2+B^2(z))}\quad\quad ({\rm by\ letting}\ \lambda=A(z)\eta)\nonumber\\
&=&\int_{\mathbf{R}^n}\frac{A(z)}{B^{2}(z)}|\hat\nu(z)|^2dz\int_{\frac{\Lambda}{A(z)}\cap\{1\le\eta\}}\frac{(\frac{B(z)}{A(z)})^2d\eta}{f(A(z)\eta)
(\eta^2+(\frac{B(z)}{A(z)})^2)}\nonumber\\
&\ge&\int_{\mathbf{R}^n}\frac{A(z)}{B^{2}(z)}|\hat\nu(z)|^2dz\sum_{k=1}^{\infty}
\int_{\{1\vee\frac{N_{x_{k}}}{A(z)}\le \eta<
\frac{N_{x_{k}+1}}{A(z)}\}}\frac{(\frac{B(z)}{A(z)})^2d\eta}{\log(N_{x_{k}+1})
(\eta^2+(\frac{B(z)}{A(z)})^2)}\nonumber\\
&=&\int_{\mathbf{R}^n}\frac{1}{B(z)}|\hat\nu(z)|^2dz\sum_{k=1}^{\infty}
\frac{1}{\log(N_{x_{k}+1})}\int_{\{1\vee\frac{N_{x_{k}}}{A(z)}\le \eta<
\frac{N_{x_{k}+1}}{A(z)}\}}\frac{\frac{A(z)}{B(z)}d\eta}{
1+(\frac{A(z)}{B(z)}\eta)^2}\nonumber\\
&=&\int_{\mathbf{R}^n}\frac{1}{B(z)}|\hat\nu(z)|^2\sum_{k=1}^{\infty}\frac{\arctan\frac{N_{x_{k}+1}}{B(z)}-
\arctan\frac{N_{{x_{k}}}\vee A(z)}{B(z)}}{\log(N_{x_{k}+1})}dz\nonumber\\
&\ge&\sum_{k=1}^{\infty}\int_{\{z\in F:\,N_{x_{k}}\le B(z)<
N_{x_{k}+1}\}}\frac{1}{B(z)}|\hat\nu(z)|^2\frac{\arctan\frac{N_{x_{k}+1}}{B(z)}-
\arctan(\frac{N_{x_{k}}}{B(z)}\vee\frac{1}{2})}{\log(N_{x_{k}+1})}dz\nonumber\\
&=&\sum_{k=1}^{\infty}\left\{\int_{\{z\in F:\,N_{x_{k}}\le B(z)<
\frac{N_{x_{k}+1}}{2}\}}\frac{1}{B(z)}|\hat\nu(z)|^2\frac{\arctan\frac{N_{x_{k}+1}}{B(z)}-
\arctan(\frac{N_{x_{k}}}{B(z)}\vee\frac{1}{2})}{\log(N_{x_{k}+1})}dz\right.\nonumber\\
&&+\left.\int_{\{z\in F:\,\frac{N_{x_{k}+1}}{2}\le B(z)<
N_{x_{k}+1}\}}\frac{1}{B(z)}|\hat\nu(z)|^2\frac{\arctan\frac{N_{x_{k}+1}}{B(z)}-
\arctan(\frac{N_{x_{k}}}{B(z)}\vee\frac{1}{2})}{\log(N_{x_{k}+1})}dz\right\}\nonumber\\
&\ge&\sum_{k=1}^{\infty}\left\{\int_{\{z\in F:\,N_{x_{k}}\le B(z)<
\frac{N_{x_{k}+1}}{2}\}}\frac{1}{B(z)}|\hat\nu(z)|^2\frac{\arctan2-
\arctan1}{\log(N_{x_{k}+1})}\right.dz\nonumber\\
&&+\left.\int_{\{z\in F:\,\frac{N_{x_{k}+1}}{2}\le B(z)<
N_{x_{k}+1}\}}\frac{1}{B(z)}|\hat\nu(z)|^2\frac{\arctan1-
\arctan\frac{1}{2}}{\log(N_{x_{k}+1})}dz\right\}\nonumber\\
&\ge&\sum_{k=1}^{\infty}\int_{\{z\in F:\,N_{x_{k}}\le B(z)<
{N_{x_{k}+1}}\}}\frac{1}{B(z)}|\hat\nu(z)|^2\frac{\arctan1-
\arctan\frac{1}{2}}{\log(N_{x_{k}+1})}dz\nonumber\\
&\ge&\frac{\arctan1- \arctan\frac{1}{2}}{\varsigma}
\sum_{k=1}^{\infty}\int_{\{z\in F:\,y_k\le B(z)<(y_k)^{\varsigma}
\}}\frac{1}{B(z)\log(B(z))}|\hat \nu(z)|^2dz.\nonumber
\end{eqnarray}
Therefore, (\ref{proof-them2.3-ii-a-1}) holds and the proof of (ii) is complete.\hfill\fbox

\subsection{Proof of Theorem \ref{thm2.4}}

The proof of Theorem \ref{thm2.4} is similar to that of Theorem \ref{thm2.3} but is more delicate. In the proof below, we will concentrate on the differences.

(i) Let $\nu$ be an arbitrary finite measure on $\mathbf{R}^n$ of
finite 1-energy. We choose a constant $\varsigma$ and a sequence
$\{x_k\}$ described as in Condition ${\rm (C^{\log\log})}$. We
assume without loss of generality that $x_1>2$. We will show that
(\ref{proof-them2.3-b}) holds.

For $\lambda>1$, we define
$$
g(\lambda):=\sum_{k=1}^{\infty}\frac{1_{\{N_{x_k+\frac{1}{2}}\le\lambda<
N_{x_k+\frac{3}{4}}\}}}{x_k\log\lambda}.
$$
Similar to (\ref{proof-them2.3-c}), we can show that
\begin{eqnarray}\label{wei513}
\int_{N_{x_1}}^{\infty}\frac{g(\lambda)}{\lambda}d\lambda
&=&\sum_{k=1}^{\infty}\frac{1}{x_k}\int_{N_{x_k+\frac{1}{2}}}^{N_{x_k+\frac{3}{4}}}\frac{1}{\lambda\log \lambda}d\lambda\nonumber\\
&=&\frac{1}{4}\log\varsigma\sum_{k=1}^{\infty}\frac{1}{x_k}\nonumber\\
&=&\infty.
\end{eqnarray}
Further, similar to (\ref{proof-them2.3-d}), we obtain that
\begin{eqnarray}\label{wei511}
& &\int_{N_{x_1}}^{\infty}\frac{g(\lambda)d\lambda}{\lambda
}\int_{\{B(z)\ge
N_{x_1}\}}\frac{\lambda}{\lambda^2+B^2(z)}|\hat\nu(z)|^2dz\nonumber\\
&\le&\sum_{l=1}^{\infty}\int_{\{N_{x_{l}}\le B(z)<
N_{(x_{l})+1}\}}\frac{1}{B(z)}|\hat\nu(z)|^2\left\{\frac{N_{x_{(l-1)}
+\frac{3}{4}}}{B(z)}\sum_{k=1}^{l-1}\frac{1}{x_k\log(N_{x_k+\frac{1}{2}})}\right.\nonumber\\
& &\ \ \ \
\left.+\frac{\pi}{2x_l\log(N_{x_l+\frac{1}{2}})}+\sum_{k=l+1}^{\infty}
\frac{\frac{B(z)}{N_{x_k+\frac{1}{2}}}}{x_k\log(N_{x_k+\frac{1}{2}})}\right\}dz\nonumber\\
&&+\sum_{l=1}^{\infty}\int_{\{N_{(x_l)+1}\le B(z)<
N_{x_{(l+1)}}\}}\frac{1}{B(z)}|\hat\nu(z)|^2\left\{\frac{N_{(x_{l})+\frac{3}{4}}}{B(z)}
\sum_{k=1}^{l}\frac{1}{x_k\log(N_{x_k+\frac{1}{2}})}\right.\nonumber\\
& &\ \ \ \
\left.+\sum_{k=l+1}^{\infty}\frac{\frac{B(z)}{N_{x_k+\frac{1}{2}}}}{x_k\log(N_{x_k+\frac{1}{2}})}\right\}dz.\
\ \ \ \ \ \ \
\end{eqnarray}

When $N_{x_{l}}\le B(z)<N_{(x_l)+1}$, by the definition of $N_x$ and the assumption that $x_1>2$, we get
\begin{eqnarray}\label{wei512}
\frac{1}{\log\log(B(z))}
&>&\frac{1}{\log\log(N_{(x_l)+1})}
\nonumber\\
&=&\frac{1}{\log\log(\varsigma^{(\varsigma^{(x_l)+1})})}\nonumber\\
&=&\frac{1}{(x_l+1)\log\varsigma+\log\log\varsigma}\nonumber\\
&>&\frac{1}{(x_l+2)\log\varsigma}\nonumber\\
&>&\frac{1}{2(\log\varsigma)x_l}.
\end{eqnarray}

We fix a $\delta>0$. Then, we obtain by (\ref{wei511}), (\ref{wei512}), (\ref{proof-them2.3-e})-(\ref{proof-them2.3-k}), the fact that $x_k>2$ for $k\ge 1$, (\ref{new345}) and Lemma \ref{slight2} that
\begin{eqnarray}\label{wei601}
& &\int_{N_{x_1}}^{\infty}\frac{g(\lambda)d\lambda}{\lambda
}\int_{\{B(z)\ge
N_{x_1}\}}\frac{\lambda}{\lambda^2+B^2(z)}|\hat\nu(z)|^2dz\nonumber\\
&\le&\frac{1}{2\varsigma^{x_1+\frac{1}{2}}(\log\varsigma)(1-\frac{1}{\varsigma})}
\int_{\mathbf{R}^n}\frac{1}{B(z)}|\hat\nu(z)|^2
\left(\frac{1}{(B(z))^{1-\varsigma^{-\frac{1}{4}}}}+
\frac{1}{(B(z))^{\varsigma^{\frac{1}{2}}-1}}\right)dz\nonumber\\
&&+{\pi(\log \varsigma)
\varsigma^{\frac{1}{2}}}\sum_{l=1}^{\infty}\int_{\{N_{x_{l}}\le
B(z)< N_{(x_{l})+1}\}}\frac{1}{B(z)\log(B(z))[\log\log(B(z))]}|\hat\nu(z)|^2dz\nonumber\\
&\le&D_1\int_{\mathbf{R}^n}
\frac{1}{B(z)\log(2+B(z))[\log\log(2+B(z))]^{1+\delta}}|\hat
\nu(z)|^2dz\nonumber\\
&&+ {\pi(\log \varsigma)
\varsigma^{\frac{1}{2}}}\sum_{l=1}^{\infty}\int_{\{y_l\le
B(z)< (y_l)^{\varsigma}\}}\frac{1}{B(z)\log(B(z))[\log\log(B(z))]}|\hat\nu(z)|^2dz\nonumber\\
&<&\infty,
\end{eqnarray}
where $D_1$ is a positive constant depending only on
$\varsigma$.

By (\ref{wei513}) and (\ref{wei601}), we obtain
(\ref{proof-them2.3-b}). Therefore, the proof of (i) is complete.
\vskip 0.2cm (ii) Suppose that $X$ satisfies (H). Let $\nu$ be an
arbitrary finite measure on $\mathbf{R}^n$ of finite 1-energy and
$\varsigma>1$ be a constant.  We choose a sequence of increasing
natural numbers $\{p_k\}$ satisfying
\begin{equation}\label{dfgo}
c(\lambda)\le\frac{1}{2^k},\ \ {\rm if}\ \lambda\ge p_k.
\end{equation}
For $m=1,2\dots$, we set $x_{m,1}=p_m+2$ and choose $k_m$ such
that
\begin{equation}\label{dec}
1\le\frac{1}{p_m+2}+\frac{1}{p_m+4}+\cdots+\frac{1}{p_m+2k_m}\le
2.
\end{equation}
Define $x_{m,l_m}=p_m+2l_m$ for $1\le l_m\le k_m$. We require
without loss of generality that
\begin{equation}\label{qguai}
x_{1,1}>1-\frac{2\log\log\varsigma}{\log\varsigma},
\end{equation}
and
$$
N^{\varsigma}_{x_{1,1}}>e\ \ {\rm and}\ \ x_{m,k_m}<p_{m+1},\ \
m=1,2,\dots
$$
Denote
$$ F:=\{z\in\mathbf{R}^n:B(z)\ge 2A(z)\}.
$$
We will show below that
$$ \sum_{m=1}^{\infty}\sum_{l_m=1}^{k_m}\int_{\{z\in F:N_{x_{m,l_m}}\le B(z)< N_{x_{m,l_m}+1}\}}
\frac{1}{B(z)\log(B(z))[\log\log(B(z))]}|\hat \nu(z)|^2dz<\infty.
$$

We define
$$
f(\lambda)=k(\log(N_k)-\log(N_{k-1})),\ \ {\rm when}\ N_{k-1}\le
\lambda<N_k,\ k\ge 2.
$$
Set $$
\Lambda:=\bigcup_{m=1}^{\infty}\bigcup_{l_m=1}^{k_m}\{\lambda:N_{x_{m,l_m}}\le
\lambda<N_{x_{m,l_m}+1}\}.
$$
Similar to (\ref{proof-them2.3-ii-b}), we can show that
$$\int_{\Lambda}\frac{c(\lambda)}{\lambda
f(\lambda)}d\lambda<\infty.
$$
Further, we obtain by (\ref{dfgo}) and (\ref{dec}) that
\begin{eqnarray}\label{wsd}
\infty&>&\int_{\Lambda}\frac{d\lambda}{\lambda
f(\lambda)}\int_{\mathbf{R}^n}\frac{\lambda}{\lambda^2+B^2(z)}|\hat\nu(z)|^21_{\{A(z)\le\lambda\}}dz\nonumber\\
&=&\int_{\mathbf{R}^n}|\hat\nu(z)|^2dz\int_{\Lambda\cap\{A(z)\le\lambda\}}\frac{d\lambda}{f(\lambda)
(\lambda^2+B^2(z))}\nonumber\\
&=&\int_{\mathbf{R}^n}|\hat\nu(z)|^2dz\int_{\frac{\Lambda}{A(z)}\cap\{1\le\eta\}}\frac{A(z)d\eta}{f(A(z)\eta)
(A^2(z)\eta^2+B^2(z))}\nonumber\\
&=&\int_{\mathbf{R}^n}\frac{A(z)}{B^{2}(z)}|\hat\nu(z)|^2dz\int_{\frac{\Lambda}{A(z)}\cap\{1\le\eta\}}\frac{(\frac{B(z)}{A(z)})^2d\eta}{f(A(z)\eta)
(\eta^2+(\frac{B(z)}{A(z)})^2)}\nonumber\\
&=&\int_{\mathbf{R}^n}\frac{A(z)}{B^{2}(z)}|\hat\nu(z)|^2dz\sum_{m=1}^{\infty}\sum_{l_m=1}^{k_m}
\int_{\{1\vee\frac{N_{x_{m,l_m}}}{A(z)}\le \eta<
\frac{N_{x_{m,l_m}+1}}{A(z)}\}}\nonumber\\
& &\ \ \ \ \
\cdot\frac{(\frac{B(z)}{A(z)})^2d\eta}{({x_{m,l_m}+1})(\log(N_{x_{m,l_m}+1})-\log(N_{{x_{m,l_m}}}))
(\eta^2+(\frac{B(z)}{A(z)})^2)}\nonumber\\
&=&\int_{\mathbf{R}^n}\frac{1}{B(z)}|\hat\nu(z)|^2dz\sum_{m=1}^{\infty}\sum_{l_m=1}^{k_m}\frac{\arctan\frac{N_{x_{m,l_m}+1}}{B(z)}-
\arctan\frac{N_{{x_{m,l_m}}}\vee A(z)}{B(z)}}{({x_{m,l_m}+1})(\log(N_{x_{m,l_m}+1})-\log(N_{{x_{m,l_m}}}))}\nonumber\\
&\ge&\sum_{m=1}^{\infty}\sum_{l_m=1}^{k_m}\int_{\{z\in
F:\,N_{x_{m,l_m}}\le B(z)<
N_{x_{m,l_m}+1}\}}\frac{1}{B(z)}|\hat\nu(z)|^2\nonumber\\
& &\ \ \ \ \ \cdot\frac{\arctan\frac{N_{x_{m,l_m}+1}}{B(z)}-
\arctan\frac{N_{x_{m,l_m}}}{B(z)}\vee\frac{1}{2})}{({x_{m,l_m}+1})(\log(N_{x_{m,l_m}+1})-\log(N_{{x_{m,l_m}}}))}dz\nonumber\\
&=&\sum_{m=1}^{\infty}\sum_{l_m=1}^{k_m}\left\{\int_{\{z\in
F:\,N_{x_{m,l_m}}\le B(z)<
\frac{N_{x_{m,l_m}+1}}{2}\}}\frac{1}{B(z)}|\hat\nu(z)|^2\right.\nonumber\\
& &\ \ \ \ \ \cdot\left.\frac{\arctan\frac{N_{x_{m,l_m}+1}}{B(z)}-
\arctan(\frac{N_{x_{m,l_m}}}{B(z)}\vee\frac{1}{2})}{({x_{m,l_m}+1})(\log(N_{x_{m,l_m}+1})-\log(N_{{x_{m,l_m}}}))}dz\right.\nonumber\\
&&+\left.\int_{\{z\in F:\,\frac{N_{x_{m,l_m}+1}}{2}\le B(z)<
N_{x_{m,l_m}+1}\}}\frac{1}{B(z)}|\hat\nu(z)|^2\right.\nonumber\\
& &\ \ \ \ \ \cdot\left.\frac{\arctan\frac{N_{x_{m,l_m}+1}}{B(z)}-
\arctan(\frac{N_{x_{m,l_m}}}{B(z)}\vee\frac{1}{2})}{({x_{m,l_m}+1})(\log(N_{x_{m,l_m}+1})-\log(N_{{x_{m,l_m}}}))}dz\right\}\nonumber\\
&\ge&\sum_{m=1}^{\infty}\sum_{l_m=1}^{k_m}\left\{\int_{\{z\in
F:\,N_{x_{m,l_m}}\le B(z)<
\frac{N_{x_{m,l_m}+1}}{2}\}}\frac{1}{B(z)}|\hat\nu(z)|^2\right.\nonumber\\
& &\ \ \ \ \ \cdot\left.\frac{\arctan2-
\arctan1}{({x_{m,l_m}+1})(\log(N_{x_{m,l_m}+1})-\log(N_{{x_{m,l_m}}}))}dz\right.\nonumber\\
&&+\left.\int_{\{z\in F:\,\frac{N_{x_{m,l_m}+1}}{2}\le B(z)<
N_{x_{m,l_m}+1}\}}\frac{1}{B(z)}|\hat\nu(z)|^2\right.\nonumber\\
& &\ \ \ \ \ \cdot\left.\frac{\arctan1-
\arctan\frac{1}{2}}{({x_{m,l_m}+1})(\log(N_{x_{m,l_m}+1})-\log(N_{{x_{m,l_m}}}))}dz\right\}\nonumber\\
&\ge&\sum_{m=1}^{\infty}\sum_{l_m=1}^{k_m}\int_{\{z\in
F:\,N_{x_{m,l_m}}\le B(z)<
{N_{x_{m,l_m}+1}}\}}\frac{1}{B(z)}|\hat\nu(z)|^2\nonumber\\
& &\ \ \ \ \ \cdot\frac{\arctan1-
\arctan\frac{1}{2}}{({x_{m,l_m}+1})(\log(N_{x_{m,l_m}+1})-\log(N_{{x_{m,l_m}}}))}dz\nonumber\\
&\ge&\frac{(\log \varsigma)(\arctan1-
\arctan\frac{1}{2})}{2\varsigma}
\sum_{m=1}^{\infty}\sum_{l_m=1}^{k_m}\int_{\{z\in
F:\,N_{x_{m,l_m}}\le B(z)< N_{x_{m,l_m}+1}\}}\nonumber\\
& &\ \ \ \ \ \cdot \frac{1}{B(z)\log(B(z))[\log\log(B(z))]}|\hat
\nu(z)|^2dz,
\end{eqnarray}
where (\ref{qguai}) has been used to obtain the last inequality.

Set $\{x_k\}=\{x_{m,l_m}:1\le l_m\le k_m\}$. Then, we obtain by
(\ref{wsd}) that
$$
\sum_{k=1}^{\infty}\int_{\{z\in F:\,N_{x_k}\le B(z)< N_{x_k+1}\}}
\frac{1}{B(z)\log(B(z))[\log\log(B(z))]}|\hat \nu(z)|^2dz<\infty.
$$
Therefore, the proof of (ii) is complete by (\ref{bn111}). \hfill\fbox

%%%%%%%%%%%%%%%%%%%%%%%%%%%%%%%%%%%%%%%%%%%%%%%%%%%%%%%%%%%%%%%%%%%%%%%
\section{New examples of L\'evy processes satisfying (H)}

As applications of Theorems \ref{thm2.3} and \ref{thm2.4}, we will
present in this section some new examples of L\'evy processes
satisfying Hunt's hypothesis (H).

\subsection{Application of Theorem \ref{thm2.3}}

Theorem \ref{thm2.3} provides a new necessary and sufficient
condition for  the validity of  (H) for L\'evy processes.
Different from the classical Kanda-Forst condition (\ref{Kanda})
and Rao's condition (\ref{Rao}), our Condition (${\rm C^{\log}}$)
only requires  that Im$(\psi)$ is partially well-controlled by
$1+{\rm Re}(\psi)$. This weaker condition is fulfilled by more
general L\'evy processes and reveals the more essential reason for
the validity of (H).

%In this section, we will present some new examples of L\'evy processes
%satisfying (H).

First, we give the following consequence of
Theorem \ref{thm2.3}.

%The purpose of these preliminary examples is to exhibit the
%strength of the criteria obtained in Section 2. Further examples
%should be expected. In fact, based on Theorems \ref{thm2.3} and
%\ref{add1}, there are two ways to further explore Getoor's
%conjecture. One the one hand, we can construct more general
%examples of L\'evy processes satisfying Condition (${\rm
%C^{\log}}$), or Condition (${\rm C^{\log\log}}$), and hence
%satisfy (H). One the other hand, we can consider if there exists a
%(maybe pathological) L\'evy process which does not satisfy
%Condition (${\rm C^{\log}}$) and hence does not satisfy (H).

\begin{prop}\label{partial}
$X$ satisfies (H)  if the following conditions hold:

(i) $ \frac{1}{c}|z|^{\alpha}\le A(z)\le B(z)\le c|z|^{\beta}$ for
$|z|\ge 1$, where $0<\alpha<\beta\le 2$ and $c>1$ are constants.

(ii) There\ exist\ $\varsigma>1$, $\kappa>0$, and\ a\ sequence\
$\{z_k\}$ such\ that $z_1>1$,
$c^{\frac{\varsigma+1}{\alpha}}z_k^{\frac{\varsigma\beta}{\alpha}}<z_{k+1}$,
$k\in \mathbf{N}$, and
$$
B(z)\le \kappa A(z)\log(B(z)),\ \ {\rm for}\ z_k\le |z|<
c^{\frac{\varsigma+1}{\alpha}}z_k^{\frac{\varsigma\beta}{\alpha}},\
k\in\mathbf{N}.
$$
\end{prop}

\noindent \emph{Proof.} Suppose that conditions (i) and (ii) hold.  By (i) and Hartman and
Wintner \cite{Hart}, we know that
$X$ has bounded continuous transition densities.

Let $\nu$ be a finite measure  on
$\mathbf{R}^n$ of finite 1-energy. Set \begin{equation}\label{vbn1} y_k=cz_k^{\beta},\ \
k\in\mathbf{N}.
\end{equation}
Then, $y_1>c>1$ and
$$
z_{k+1}>cz_k^{\varsigma}>cz_k,\ \ k\in \mathbf{N}.
$$
Hence $z_k\uparrow\infty$ and thus $y_k\uparrow\infty$ as $k\to\infty$.

By (i) and (\ref{vbn1}), we find that when $y_k\leq B(z)<(y_k)^{\varsigma}$,
\begin{eqnarray}\label{partial-a}
z_k\le |z|<
c^{\frac{\varsigma+1}{\alpha}}z_k^{\frac{\varsigma\beta}{\alpha}}.
\end{eqnarray}
Then, we obtain by (\ref{partial-a}),  $y_1>c$, conditions (i) and (ii) that
\begin{eqnarray*}
& &\sum_{k=1}^{\infty}\int_{\{y_k\le B(z)< (y_k)^{\varsigma}\}}
\frac{1}{B(z)\log(B(z))}|\hat \nu(z)|^2dz\\
&=&\sum_{k=1}^{\infty}\int_{\{y_k\le B(z)< (y_k)^{\varsigma},|z|\geq 1\}}
\frac{1}{B(z)\log(B(z))}|\hat \nu(z)|^2dz\\
&&+\sum_{k=1}^{\infty}\int_{\{y_k\le B(z)< (y_k)^{\varsigma},|z|< 1\}}
\frac{1}{B(z)\log(B(z))}|\hat \nu(z)|^2dz\\
&\le&
\sum_{k=1}^{\infty}\int_{\{z_k\le |z|<
c^{\frac{\varsigma+1}{\alpha}}z_k^{\frac{\varsigma\beta}{\alpha}}\}}
\frac{\kappa A(z)}{B^{2}(z)}|\hat
\nu(z)|^2dz+\int_{\{z:\,|z|<1,\,B(z)\ge c\}}\frac{1}{B(z)\log(B(z))}|\hat \nu(z)|^2dz\\
&\le&\int_{\mathbf{R}^n}\frac{\kappa A(z)}{B^{2}(z)}|\hat
\nu(z)|^2dz+\frac{(\nu(\mathbf{R}^n))^2}{c\log c}\int_{\{|z|<1\}}dz\\
 &<&\infty.
\end{eqnarray*}
Therefore, Condition (${\rm C^{\log}}$) holds and the proof is
complete by Theorem \ref{thm2.3}(i).\hfill\fbox

\begin{rem}
Blumenthal and Getoor introduced in \cite{BG61} different indices
for L\'evy processes on $\mathbf{R}^n$. In particular, they
defined
$$
\beta=\inf\left\{\alpha\ge 0:
\frac{|\psi(z)|}{|z|^{\alpha}}\rightarrow 0\ {\rm as}\
|z|\rightarrow\infty\right\},
$$
and
$$
\beta''=\sup\left\{\alpha\ge 0: \frac{{\rm
Re}\psi(z)}{|z|^{\alpha}}\rightarrow \infty\ {\rm as}\
|z|\rightarrow\infty\right\}.
$$
Proposition \ref{partial} provides a sufficient condition for the
validity of (H) in the case that $\beta''<\beta$.

Note that in condition (ii) of Proposition \ref{partial}, $z_{k+1}$
can be chosen to be much bigger than $z_k$ for each
$k\in\mathbf{N}$. Hence Proposition \ref{partial} can be used to
construct a class of L\'evy processes with indices $\beta''<\beta$
and satisfying (H). As a concrete example, we will give in Example
\ref{example1} below a class of subordinators satisfying (H).
Recall that by virtue of the following remarkable result of Glover
and Rao, whenever we find a new class of subordinators
satisfying (H), we obtain a new class of time-changed Markov
processes satisfying (H).

\begin{prop}\label{GR} (Glover and Rao \cite{GR86}) Let $(X_t)_{t\ge 0}$ be a standard Markov process on a locally compact space with a
countable base and $(T_t)_{t\ge 0}$ be an independent subordinator
satisfying Hunt's hypothesis (H). Then $(X_{T_t})_{t\ge 0}$
satisfies (H).
\end{prop}

\end{rem}

\begin{exa}\label{example1}
Let $0<\alpha<\beta<1$, $c_1>1$, $\varsigma>1$ and $0<\kappa_1\le
c_1$. Denote $$
c:=c_1\left(8+\frac{1}{1-\beta}+\frac{2}{\beta}\right).
$$
We choose a sequence $\{z_k\}$ satisfying $z_1>1$,
$c^{\frac{\varsigma+1}{\alpha}}z_k^{\frac{\varsigma\beta}{\alpha}}<z_{k+1}$,
$k\in \mathbf{N}$.

Let $X$ be a pure jump subordinator with L\'evy measure
$\mu(dx):=\rho(x)dx$. Suppose $\rho$ satisfies the following
conditions:

(i) $\frac{1}{c_1x^{1+\alpha}}\le\rho(x)\le
\frac{c_1}{x^{1+\beta}}$ for $0<x\le 1$; and $\rho(x)=0$,
otherwise.

(ii) $\rho(x)\ge \frac{\kappa_1}{x^{1+\beta}}$ for
$\frac{1}{2c^{\frac{\varsigma+1}{\alpha}}z_k^{\frac{\varsigma\beta}{\alpha}}}\le
x< \frac{1}{z_{k}}$, $k\in \mathbf{N}$.

\noindent Then $X$ satisfies (H).

In fact, for $z\in \mathbf{R}$ with $|z|\ge1$, we have
\begin{eqnarray}\label{exam-3.4-a}
A(z)&=&1+\int_0^{\infty}(1-\cos(zx))\rho(x)dx\nonumber\\
&\ge&\int_{0}^{\frac{1}{|z|}}(1-\cos(zx))
\frac{1}{c_1x^{1+\alpha}}dx\nonumber\\
&\ge&\int_{0}^{\frac{1}{|z|}}
\frac{|z|^2x^2}{4c_1x^{1+\alpha}}dx\nonumber\\
&=&\frac{1}{4c_1(2-\alpha)}|z|^{\alpha}\nonumber\\
&>&\frac{1}{8c_1}|z|^{\alpha}\nonumber\\
&>&\frac{1}{c}|z|^{\alpha},
\end{eqnarray}
and
\begin{eqnarray}\label{exam-3.4-b}
B(z)&=&|1+\psi(z)|\nonumber\\
&\leq& 1+{\rm Re}\psi(z)+\left|{\rm Im}\psi(z)\right|\nonumber\\
&=&1+\int_0^1(1-\cos(zx))\rho(x)dx+\left|\int_0^1\sin(zx)\rho(x)dx\right|\nonumber\\
&\leq &1+\int_0^1(1-\cos(zx))\frac{c_1}{x^{1+\beta}}dx+\int_0^1|\sin(zx)|\frac{c_1}{x^{1+\beta}}dx\nonumber\\
&\le&1+\int_{0}^{\frac{1}{|z|}}
\left(\frac{|z|^2x^2}{2}+|z|x\right)\frac{c_1}{x^{1+\beta}}dx+2\int_{\frac{1}{|z|}}^{1}\frac{c_1}{x^{1+\beta}}dx\nonumber\\
&=&1+\frac{c_1}{2(2-\beta)}|z|^{\beta}+\frac{c_1}{1-\beta}|z|^{\beta}+\frac{2c_1}{\beta}(|z|^{\beta}-1)\nonumber\\
&=&c_1\left(\frac{1}{2(2-\beta)}+\frac{1}{1-\beta}+\frac{2}{\beta}\right)|z|^{\beta}+\left(1-\frac{2c_1}{\beta}\right)\nonumber\\
&<&c|z|^{\beta}.
\end{eqnarray}

For $z_k\le
|z|<c^{\frac{\varsigma+1}{\alpha}}z_k^{\frac{\varsigma\beta}{\alpha}}$,
$k\in \mathbf{N}$, by condition (ii) and (\ref{exam-3.4-b}), we get
\begin{eqnarray*}
A(z) &\ge&\int_{\frac{1}{2|z|}}^{\frac{1}{|z|}}(1-\cos(zx))
\frac{\kappa_1}{x^{1+\beta}}dx\\
&\ge&\int_{\frac{1}{2|z|}}^{\frac{1}{|z|}}
\frac{\kappa_1|z|^2x^2}{4x^{1+\beta}}dx\\
&=&\frac{\kappa_1}{4(2-\beta)}\left(1-\frac{1}{2^{2-\beta}}\right)|z|^{\beta}\\
&\ge&\frac{\kappa_1}{16}|z|^{\beta}\\
&\geq& \frac{\kappa_1}{16c}B(z).
\end{eqnarray*}

By (\ref{exam-3.4-a}) and Hartman and
Wintner \cite{Hart}, $X$ has bounded continuous transition
densities.  Hence conditions (i) and (ii) of Proposition \ref{partial} are
fulfilled and therefore $X$ satisfies (H).
\end{exa}

\subsection{Application of Theorem \ref{thm2.4}}

We define
\begin{eqnarray*}
& &{\rm  Condition\ (C^0)}:\ {\rm For\ any\ finite\
measure}\ \nu\ {\rm on}\ \mathbf{R}^n\ {\rm of\ finite\ 1{\textrm{-}}energy,}\\
& &\ \ \ \ \ \int_{\mathbf{R}^n} \frac{1}{B(z)\log
(2+B(z))[\log\log(2+B(z))]}|\hat \nu(z)|^2dz<\infty,
\end{eqnarray*}
and
\begin{eqnarray*}
& &{\rm Condition\ (C^{B/A}):}\ {\rm There\ exists\ a\ constant}\ C>0\ {\rm such\ that}\\
& &\ \ \ \ \ B(z)\le CA(z)\log (2+B(z))[\log\log(2+B(z))],\ \
\forall z\in \mathbf{R}^n.
\end{eqnarray*}

As a direct consequence of Theorem \ref{thm2.4}, we obtain the
following new sufficient conditions for the validity of Hunt's
hypothesis (H).
 \vskip 0.5cm
 \begin{prop}\label{slight} Condition ($C^{B/A}$)
$\Rightarrow$ Condition ($C^0$) $\Rightarrow$ (H).
\end{prop}

 \vskip 0.2cm By
Proposition \ref{slight}, we get the following  result.

\begin{cor}\label{cor-2.6}
Let $X$ be a L\'evy process on $\mathbf{R}$. Suppose that
\begin{equation}\label{cor-2.6-a}
\liminf_{|z|\to
\infty}\frac{|\psi(z)|}{|z|(\log\log|z|)^{\delta}}>0
\end{equation}
for some constant $\delta>0$. Then $X$ satisfies (H).
\end{cor}

\noindent \emph{Proof.} By (\ref{cor-2.6-a}), there exist constants $K>e$ and $\kappa>0$ such that
\begin{eqnarray*}
\frac{|\psi(z)|}{|z|(\log\log|z|)^{\delta}}\geq \kappa,\ \ {\rm if}\ |z|\ge K.
\end{eqnarray*}
Hence, when $|z|\geq K$, we have
$$B(z)=|1+\psi(z)|\geq |\psi(z)|\geq \kappa |z|(\log\log|z|)^{\delta}.
$$

Define $G=\max\{K,e^e,\frac{2e}{\kappa}\}$ and
$$
C_1=\int_{|z|<G} \frac{1}{B(z)\log(2+B(z))[\log\log(2+B(z))]}dz.
$$
Then, we have
\begin{eqnarray*}
&&\int_{\mathbf{R}^n} \frac{1}{B(z)\log(2+B(z))[\log\log(2+B(z))]}dz\\
&=&C_1+\int_{|z|\geq G} \frac{1}{B(z)\log(2+B(z))[\log\log(2+B(z))]}dz\\
&\le& C_1+\int_{|z|\geq G} \frac{1}{\kappa |z|(\log\log|z|)^{\delta}\log(2+\kappa |z|(\log\log|z|)^{\delta})[\log\log(2+\kappa |z|(\log\log|z|)^{\delta})]}dz\\
&\le& C_1+\int_{|z|\geq G} \frac{1}{\kappa |z|(\log\log|z|)^{\delta}\log(\kappa |z|)[\log\log(\kappa |z|)]}dz\\
&\le& C_1+C_2\int_{|z|\geq G} \frac{1}{|z|\log|z|(\log\log|z|)^{1+\delta}}dz\\
&<&\infty,
\end{eqnarray*}
where $C_2$ is a positive constant depending on $K$ and $\kappa$. Therefore, Condition $(C^0)$ holds and the proof is complete by Proposition \ref{slight}.\hfill\fbox

 \vskip 0.5cm
\begin{exa}\label{exm3.5}
Let $X$ be a L\'{e}vy process  on $\mathbf{R}$ with
L\'{e}vy-Khintchine exponent $(a,Q,\mu)$. Suppose that there exist
constants $\delta>0$ and $c>0$ such that
$$
d\mu\geq \frac{c[\log(-\log|x|)]^{\delta}}{x^2}dx \ \ { on}\ \
\left\{x\in \mathbf{R}: 0<|x|<\frac{1}{e}\right\}.$$
Then $X$
satisfies (H).

In fact, for $|z|\geq e$, we have
\begin{eqnarray*}
{\rm Re}\psi(z)&=&\frac{1}{2}Qz^2+\int_{\mathbf{R}}(1-\cos(zx))\mu(dx)\\
&\geq&\int_{\{|x|<\frac{1}{e}\}}(1-\cos(zx))\frac{c[\log(-\log|x|)]^{\delta}}{x^2}dx\\
&\geq&2\int_0^{\frac{1}{|z|}}(1-\cos(zx))\frac{c[\log(-\log|x|)]^{\delta}}{x^2}dx\\
&\geq &2\int_0^{\frac{1}{|z|}}\frac{|z|^2x^2}{2}\cdot\frac{c[\log(-\log|x|)]^{\delta}}{x^2}dx\\
&=&c|z|^2\int_0^{\frac{1}{|z|}}[\log(-\log|x|)]^{\delta}dx\\
&\geq&c|z|^2\int_0^{\frac{1}{|z|}}(\log\log|z|)^{\delta}dx\\
&=&c|z|(\log\log|z|)^{\delta},
\end{eqnarray*}
which implies  (\ref{cor-2.6-a}). Therefore, $X$
satisfies (H) by Corollary \ref{cor-2.6}. Note that in this
example it does not matter if $a$ or $Q$ equals 0.
\end{exa}

\begin{rem} Corollary \ref{cor-2.6} and Example \ref{exm3.5} extend the corresponding results of \cite[Proposition 4.10 and Example 4.8]{HSZ15}.
\end{rem}

\section{Decomposition of subordinator into two independent ones satisfying (H)}%\setcounter{equation}{0}

In \cite{HS12}, we proved that if a subordinator satisfies (H) then it must be a pure jump
subordinator. To date, it is still unknown if all pure jump subordinators satisfy (H). We refer the interested readers to  \cite{HSZ15} for some known examples of subordinators
satisfying  (H).

\begin{defin} \label{defi-3.6}
Let $0<\alpha<\beta<1$. A pure jump subordinator $X$ is said to be of type-$(\alpha,\beta)$ if the L\'evy measure of $X$ has density, which is denoted by $\rho$, and there exists a constant $c>1$ such that
\begin{equation}\label{alphabeta}
\frac{1}{cx^{1+\alpha}}\le\rho(x)\le \frac{c}{x^{1+\beta}},\ \ \forall x\in
(0,1].\end{equation}
\end{defin}

In this section, we will apply Theorem \ref{thm2.3} to prove the following result.

\begin{thm}\label{sub1}
 Any pure jump subordinator of type-$(\alpha,\beta)$ can be decomposed into the summation of two independent pure jump subordinators of type-$(\alpha,\beta)$ such that both of them satisfy (H).
\end{thm}

\noindent \emph{Proof.}  We fix $0<\alpha<\beta<1$. Let $X$ be a pure jump subordinator of type-$(\alpha,\beta)$ such that its L\'evy density $\rho$ satisfies (\ref{alphabeta}). By \cite[Theorem 2.1 and
Corollary 2.2]{HSZ15}, big jumps have no effect on the validity of
(H). So we can assume without loss of generality that $\rho(x)=0$
for $x\in (1,\infty)$.

For $i=1,2$, we set
\begin{eqnarray}\label{2345}
\rho^i(x)=\left\{
  \begin{array}{ll}
   \frac{1}{2cx^{1+\alpha}}+\rho_i(x),&\ {\rm if}\ x\in (0,1],\\
    0,&\ {\rm if}\ x\in (1,\infty),
  \end{array}
\right.
\end{eqnarray}
where $\rho_1$ and $\rho_2$ are two non-negative functions defined
on $(0,1]$ satisfying
$$
\rho_1(x)+\rho_2(x)=\rho(x)-\frac{1}{cx^{1+\alpha}}, \ \
\forall x\in (0,1],
$$
and thus
\begin{eqnarray}\label{2346-a}
\rho^1(x)+\rho^2(x)=\rho(x), \ \ \forall x\in (0,\infty).
\end{eqnarray}
We will suitably define $\rho_1,\rho_2$ below to ensure that the pure jump subordinators $X^1$
and $X^2$ with L\'evy densities $\rho^1$ and $\rho^2$,
respectively, satisfy the requirements of the desired
decomposition. Note that by (\ref{alphabeta})-(\ref{2346-a}), for  $i=1,2$, we have
\begin{eqnarray}\label{2347}
\frac{1}{2cx^{1+\alpha}}\leq \rho^i(x)\leq \frac{c}{x^{1+\beta}},\ \  \forall x\in (0,1].
\end{eqnarray}

For $i=1,2$, we use $\psi_{X^i}$ to denote the
L\'{e}vy-Khintchine exponent of $X_i$ and use $A_{X^i}$,
$B_{X^i}$, etc. to denote the corresponding quantities of $X^i$
(cf. (\ref{ABAB})). For  $z\in \mathbf{R}$ with $|z|\ge 1$, by (\ref{2347}), $c>1$ and $0<\alpha< \beta<1$,  we obtain that
\begin{eqnarray}\label{ZA0}
A_{X^i}(z)&=&1+{\rm Re}\psi_{X^i}(z)\nonumber\\
&=&1+\int_0^{\infty}(1-\cos(zx))\rho_i(x)dx\nonumber\\
&\geq&1+\int_0^1(1-\cos(zx))\frac{1}{2cx^{1+\alpha}}dx\nonumber\\
&\ge&1+\int_{0}^{\frac{1}{|z|}}(1-\cos(zx))
\frac{1}{2cx^{1+\alpha}}dx\nonumber\\
&\ge&1+\int_{0}^{\frac{1}{|z|}}
\frac{|z|^2x^2}{8cx^{1+\alpha}}dx\nonumber\\
&=&1+\frac{1}{8c(2-\alpha)}|z|^{\alpha}\nonumber\\
&>&1+\frac{1}{16c}|z|^{\alpha},
\end{eqnarray}
and
\begin{eqnarray}\label{BB}
B_{X^i}(z)&=&|1+\psi_{X^i}(z)|\nonumber\\
&\le&1+\int_0^1(1-\cos(zx))\rho_i(x)dx+\left|\int_0^1\sin(zx)\rho_i(x)dx\right|\nonumber\\
&\leq &1+\int_0^1(1-\cos(zx))\frac{c}{x^{1+\beta}}dx+\int_0^1|\sin(zx)|\frac{c}{x^{1+\beta}}dx\nonumber\\
&\le&1+\int_{0}^{\frac{1}{|z|}}
\left(\frac{|z|^2x^2}{2}+|z|x\right)\frac{c}{x^{1+\beta}}dx+2\int_{\frac{1}{|z|}}^{1}\frac{c}{x^{1+\beta}}dx\nonumber\\
&=&1+\frac{c}{2(2-\beta)}|z|^{\beta}+\frac{c}{1-\beta}|z|^{\beta}+\frac{2c}{\beta}(|z|^{\beta}-1)\nonumber\\
&=&c\left(\frac{1}{2(2-\beta)}+\frac{1}{1-\beta}+\frac{2}{\beta}\right)|z|^{\beta}+\left(1-\frac{2c}{\beta}\right)\nonumber\\
&<&c\left(\frac{1}{2}+\frac{1}{1-\beta}+\frac{2}{\beta}\right)|z|^{\beta}.
\end{eqnarray}

Fix a constant $\varsigma>1$ and set $\varepsilon_0=1$,
$\varepsilon_1=1/2$, $\rho_1(x)=\rho(x)-\frac{1}{cx^{1+\alpha}}$,
 $\rho_2(x)=0$ for $x\in (\varepsilon_1,1]$. We
will define $\varepsilon_n,\rho_1(x),\rho_2(x)$, $x\in
(\varepsilon_n,\varepsilon_{n-1}]$ by induction. Suppose that
$\varepsilon_1,\dots,\varepsilon_n,\rho_1(x),\rho_2(x),x\in
(\varepsilon_n,1],n\ge 1$, have been defined. In the following, we will define
$\varepsilon_{n+1},\rho_1(x),\rho_2(x),x\in
(\varepsilon_{n+1},\varepsilon_n]$.

First, we
consider the case that $n$ is an odd number. Define
$d\mu^{(n)}=(\rho_1(x)1_{(\varepsilon_{n},1](x)})dx$. By the
Riemann-Lebesgue lemma,
$$
\lim_{z\rightarrow\infty}\int_{\mathbf{R}}\sin(zx)\mu^{(n)}(dx)=0.
$$
Hence there exists a positive constant, which is denoted by $z_{n+1}$, such
that
\begin{equation}\label{Lebes0}
\left|\int_{\mathbf{R}}\sin(zx)\mu^{(n)}(dx)\right|\le 1,\ \ {\rm
for}\ |z|\ge z_{n+1}.
\end{equation}
We assume without loss of generality that
\begin{equation}\label{zxc1}
z_{n+1}>\frac{1}{\varepsilon_{n}}.
\end{equation}
Define
\begin{eqnarray}
&&c_1:=c\left(\frac{1}{2}+\frac{1}{1-\beta}+\frac{2}{\beta}\right),\label{ec1}\\
&&z_{n+1}':=\left[16c(c_1z^{\beta}_{n+1})^{\zeta}\right]^{1/\alpha},\label{dcv}\\
&&\varepsilon_{n+1}:=1/(z_{n+1}')^{1/(1-\beta)}.\label{ssa}
\end{eqnarray}
For any $x\in (\varepsilon_{n+1},\varepsilon_n]$,  define
$\rho_1(x)=0$ and $\rho_2(x)=\rho(x)-\frac{1}{cx^{1+\alpha}}$. For
the case that $n$ is an even number, we can similarly define
$\varepsilon_{n+1},\rho_1(x),\rho_2(x),x\in
(\varepsilon_{n+1},\varepsilon_{n}]$ but with the places of
$\rho_1$ and $\rho_2$ switched.

Now for $i=1,2$, we let $\rho^i$ be defined as above and let $X^i$
be a pure jump subordinator with the L\'evy density $\rho^i$. By
(\ref{2346-a}), we can assume that $X^1$ and $X^2$ are independent and $X=X^1+X^2$.
We will show that $X^1$ satisfies (H). The proof of the validity of (H) for $X^2$ is
similar so we omit it.

Let $n\ge 1$ be an odd number. For $z_{n+1}\le
|z|\le z_{n+1}'$,  by (\ref{2345}), the fact that $\rho_1(x)=0$ for $x\in (\varepsilon_{n+1},\varepsilon_n]$, (\ref{2347}), (\ref{Lebes0}),
(\ref{ssa}) and (\ref{ZA0}), we obtain that
\begin{eqnarray}\label{cvbx0}
B_{X^1}(z)&=&|1+\psi_{X^1}(z)|\nonumber\\
&\le& A_{X^1}(z)+|{\rm Im}\psi_{X^1}(z)|\nonumber\\
&=&A_{X^1}(z)+\left|\int_0^1\sin(zx)\rho^1(x)dx\right|\nonumber\\
&\le &A_{X^1}(z)+\left|\int_0^{\varepsilon_{n+1}}\sin(zx)\rho_1(x)dx\right|+\left|
\int_{\varepsilon_{n}}^1\sin(zx)\rho_1(x)dx\right|\nonumber\\
& &+\left|\int_{0}^1\sin(zx)\frac{1}{2cx^{1+\alpha}}dx\right|\nonumber\\
&\le&A_{X^1}(z)+\int_0^{\varepsilon_{n+1}}\frac{c|z|x}{x^{1+\beta}}dx+1+\left(\int_{0}^{\frac{1}{|z|}}
\frac{|z|x}{2cx^{1+\alpha}}dx+\int_{\frac{1}{|z|}}^{1}\frac{1}{2cx^{1+\alpha}}dx\right)\nonumber\\
&=&A_{X^1}(z)+\frac{c|z|\varepsilon_{n+1}^{1-\beta}}{1-\beta}+1+\left(
\frac{|z|^{\alpha}}{2c(1-\alpha)}+\frac{1}{2c\alpha}(|z|^{\alpha}-1)\right)\nonumber\\
&\le& A_{X^1}(z)+\frac{c}{1-\beta}+1+\left(\frac{1}{2c(1-\alpha)}+\frac{1}{2c\alpha}\right)|z|^{\alpha}\nonumber\\
&\le&c^1A_{X^1}(z),
\end{eqnarray}
where $c^1>1$ is a constant independent of $z$ and $n$.

Let $\nu$ be a finite measure  on $\mathbf{R}$ of finite 1-energy
with respect to $X^1$. Set $y_k=c_1z_{k+1}^{\beta},\ k\in\mathbf{N}$. Then, by (\ref{zxc1})-(\ref{ssa}), we find that
$\{y_k\uparrow\infty\}$ and $y_1>c_1>1$. When $k$ is odd and $y_k\le B_{X^1}(z)<
(y_k)^{\varsigma}$, by (\ref{ZA0}), (\ref{BB}) and (\ref{dcv}),  we get
\begin{eqnarray}\label{hh-a}
z_{k+1}\leq |z|<z_{k+1}'.
\end{eqnarray}
By  (\ref{cvbx0}) and (\ref{hh-a}), we obtain  that
\begin{eqnarray*}
& &\sum_{k\ {\rm odd}}^{\infty}\int_{\{y_k\le B_{X^1}(z)<
(y_k)^{\varsigma}\}} \frac{1}{B_{X^1}(z)\log(B_{X^1}(z))}|\hat
\nu(z)|^2dz\\
&=&\sum_{k\ {\rm odd}}^{\infty}\int_{\{y_k\le B_{X^1}(z)<
(y_k)^{\varsigma},\, |z|\geq 1\}} \frac{1}{B_{X^1}(z)\log(B_{X^1}(z))}|\hat
\nu(z)|^2dz\\
&&+\sum_{k\ {\rm odd}}^{\infty}\int_{\{y_k\le B_{X^1}(z)<
(y_k)^{\varsigma},\, |z|<1\}} \frac{1}{B_{X^1}(z)\log(B_{X^1}(z))}|\hat
\nu(z)|^2dz\\\\
&\le& \sum_{k\ {\rm odd}}^{\infty}\int_{\{z_{k+1}\le
|z|<z_{k+1}'\}}\frac{1}{B_{X^1}(z)\log(c_1)}|\hat
\nu(z)|^2dz\\
& &+\int_{\{B_{X^1}(z)\ge c_1,\,|z|<1\}}\frac{1}{B_{X^1}(z)\log(B_{X^1}(z))}|\hat \nu(z)|^2dz\\
&\le&\frac{1}{\log(c_1)} \sum_{k\ {\rm odd}}^{\infty}\int_{\{z_{k+1}\le
|z|<z_{k+1}'\}} \frac{c^1A_{X^1}(z)}{B_{X^1}^{2}(z)}|\hat
\nu(z)|^2dz\nonumber\\
& &+\frac{1}{c_1\log(c_1)}\int_{\{|z|<1\}}|\hat \nu(z)|^2dz\\
&\le&\frac{c^1}{\log(c_1)}\int_{\mathbf{R}}\frac{A_{X^1}(z)}{B_{X^1}^{2}(z)}|\hat
\nu(z)|^2dz+\frac{(\nu(\mathbf{R}))^2}{c_1\log c_1}\int_{\{|z|<1\}}dz\nonumber\\
&<&\infty.
\end{eqnarray*}
Hence Condition (${\rm C^{\log}}$) is fulfilled. By (\ref{ZA0}) and Hartman and
Wintner \cite{Hart}, $X_1$ has bounded continuous transition
densities. Therefore, $X^1$ satisfies (H) by Theorem
\ref{thm2.3}(i).\hfill\fbox

 \vskip 0.5cm
\begin{rem}
By \cite[Theorem 2.1 and Corollary 2.2]{HSZ15}, Theorem \ref{sub1}
still holds if the interval $(0,1]$ in Definition \ref{defi-3.6}
is replaced with $(0,\delta]$ for any constant $\delta>0$.
\end{rem}

Theorem \ref{sub1} leads us to consider the following question:

{\it Suppose that $X_1$ and $X_2$ are two independent L\'evy processes
on $\mathbf{R}^n$ such that both of them satisfy (H). Does
$X_1+X_2$ satisfy (H)?}

We will consider this question in a forthcoming paper \cite{HS}.

\newpage\bigskip

{ \noindent {\bf\large Acknowledgments} \vskip 0.1cm  \noindent We wish to
thank the financial support of NNSFC (Grant No. 11371191) and NSERC (Grant No. 311945-2013).}

\end{document}